\documentclass[11pt,a4paper]{article}
\usepackage{amsfonts}
\usepackage{latexsym}
\usepackage{cite}
\usepackage{amsmath,amsfonts,latexsym,amssymb}
\usepackage[mathscr]{eucal}
\usepackage{cases}
\usepackage{amsthm}
\usepackage[bf,small]{caption2}
\usepackage{float,color}
\usepackage{graphicx}
\usepackage{amsmath}
\usepackage{amssymb}
\usepackage[all]{xy}
\usepackage{makecell}
\usepackage{multirow}
\usepackage{hyperref}

\newtheorem{theorem}{theorem}[section]
\newtheorem{thm}[theorem]{Theorem}
\newtheorem{lem}[theorem]{Lemma}
\newtheorem{prob}[theorem]{Problem}

\newtheorem{cor}[theorem]{Corollary}

\newtheorem{defn}[theorem]{Definition}
\newtheorem{exmp}[theorem]{Example}

\newtheorem{rmk}[theorem]{Remark}
\newtheorem{nota}[theorem]{Notation}

\begin{document}

\title{\vspace{-2cm}\textbf{On skein algebras of planar surfaces}}
\author{\Large Haimiao Chen}
\date{}
\maketitle

\begin{abstract}
  Let $R$ be a commutative ring with identity and a fixed invertible element $q^{\frac{1}{2}}$.
  Let $\mathcal{S}_n$ denote the Kauffman bracket skein algebra of the $n$-holed disk $\Sigma_{0,n+1}$ over $R$.
  When $q+q^{-1}$ is invertible, in 2000 Przytycki and Sikora found a set of $n+{n\choose 2}+{n\choose 3}$ generators for $\mathcal{S}_n$;
  we show that the ideal of defining relations among these generators is generated by relations of degree $\le6$ supported by certain subsurfaces diffeomorphic to $\Sigma_{0,k+1}$ with $k\le 6$. When $q+q^{-1}$ is not invertible, a set of $2^n-1$ generators for $\mathcal{S}_n$ was known to Bullock in 1999; we show that the ideal of defining relations is generated by relations of degree $\le 2k+2$ supported by certain subsurfaces diffeomorphic to $\Sigma_{0,k+1}$ with $k\le n$. These results are substantial progresses towards answering Problem 1.92 (J) in the Kirby's list.

  \medskip
  \noindent {\bf Keywords:} $n$-holed disk; Kauffman bracket skein algebra; quantization; presentation; defining relation   \\
  {\bf MSC2020:} 57K16, 57K31
\end{abstract}

\section{Introduction}

\begin{figure}[H]
  \centering
  \includegraphics[width=7cm]{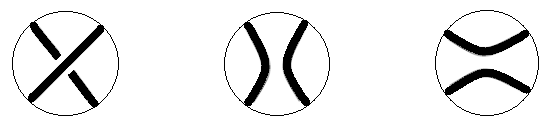}\\
  \caption{Left: $\mathbf{l}_\times$; middle: $\mathbf{l}_\infty$; right: $\mathbf{l}_0$. They are the same outside the ball.}\label{fig:local}
\end{figure}

Let $R$ be a commutative ring with identity and a fixed invertible element $q^{\frac{1}{2}}$. Given an oriented 3-manifold $M$, the {\it Kauffman bracket skein module} of $M$ over $R$, denoted by $\mathcal{S}(M;R)$, is defined as the $R$-module generated by ambient isotopy classes of (possibly empty) framed links in the interior of $M$, subject to the {\it skein relations}
\begin{align}
\mathbf{l}\sqcup\bigcirc=-(q+q^{-1})\mathbf{l}, \qquad
\mathbf{l}_\times=q^{\frac{1}{2}}\mathbf{l}_\infty+q^{-\frac{1}{2}}\mathbf{l}_0    \label{eq:local}
\end{align}
for any $\mathbf{l}$, any $0$-framed unknot $\bigcirc$ unlinked with $\mathbf{l}$, and any skein triple $(\mathbf{l}_\times,\mathbf{l}_\infty,\mathbf{l}_0)$.
By {\it skein triple} we mean a triple consisting of three links which are the same outside a small $3$-ball $B$, but inside $B$ they are as in Figure \ref{fig:local}.

As a convention, $R$ is identified with $R\emptyset\subset\mathcal{S}(M;R)$ via $\lambda\mapsto \lambda\emptyset$.

In the case $M=\Sigma\times[0,1]$ where $\Sigma$ is an oriented surface, $\mathcal{S}(M;R)$ is usually denoted as $\mathcal{S}(\Sigma;R)$, and its elements are given by $R$-linear combinations of links in $\Sigma\times(0,1)$. It is understood that each framing vector of a link is parallel to $[0,1]$ and points towards $1$. Equipped with the product defined by superposition, $\mathcal{S}(\Sigma;R)$ becomes a $R$-algebra, called the {\it Kauffman bracket skein algebra} of $\Sigma$ over $R$. Using skein relations, each element of $\mathcal{S}(\Sigma;R)$ can be written as a $R$-linear combination of multicurves, where a {\it multicurve} means a disjoint union of nonnullhomotopic simple curves and is regarded as a link in $\Sigma\times\{\frac{1}{2}\}\subset\Sigma\times(0,1)$. According to \cite[Fact 4.1]{PS00} or \cite[Corollary 4.1]{SW07},
multicurves form a free basis for the $R$-module $\mathcal{S}(\Sigma;R)$.

When $R=\mathbb{C}$ and $q^{\frac{1}{2}}=-1$, Bullock \cite{Bu97} showed that $\mathcal{S}(M;\mathbb{C})$ modulo its nilradical is isomorphic to the coordinate ring of the ${\rm SL}(2,\mathbb{C})$-character variety $\mathcal{X}(\pi_1(M))$.
In this sense, skein module is considered as a {\it quantization} of character variety.
For $\Sigma\times[0,1]$, Przytycki and Sikora \cite{PS00,PS19} showed that $\mathcal{S}(\Sigma;R)$ is a domain for any integral domain $R$;
in particular, $\mathcal{S}(\Sigma;\mathbb{C})\cong\mathbb{C}[\mathcal{X}(\pi_1(\Sigma))]$. This isomorphism was also established by Charles and March\'e \cite{CM12} through a completely different approach. Skein algebra was also proposed by Turaev \cite{Tu91} independently as a quantization of character variety.

Being fundamentally significant, the following was raised as Problem 1.92 (J) in the Kirby's list \cite{Ki97}, and also \cite[Problem 4.5]{Oh02}:
\begin{prob}[Bullock and Przytycki]\label{prob:main}
Find the structure of $\mathcal{S}(\Sigma_{g,k};\mathbb{Z}[q^{\pm\frac{1}{2}}])$.
\end{prob}
A finite set of generators was given by Bullock \cite{Bu99}, so the real problem is to determine the relations. The structure of $\mathcal{S}(\Sigma_{g,k};\mathbb{Z}[q^{\pm\frac{1}{2}}])$ for $g=0,k\le 4$ and $g=1,k\le 2$ was known to Bullock and Przytycki \cite{BP00} early in 2000. Till now it remains difficult to find all relations for general $g$ and $k$.
Recently, Cooke and Lacabanne \cite{CL25} obtained a presentation for $\mathcal{S}(\Sigma_{0,5};\mathbb{C}(q^{\frac{1}{4}}))$,
and the author \cite{Ch25} deduced a presentation for $\mathcal{S}(\Sigma_{0,5};\mathbb{Z}[q^{\pm\frac{1}{2}}])$ with essentially the same expressions as in \cite{CL25}, by a different method.
For $g\ge 1$ and $k\le 1$, Santharoubane \cite{Sa24} gave a useful criteria for showing a set of curves to generate $\mathcal{S}(\Sigma_{g,k};\mathbb{Q}(q^{\frac{1}{2}}))$, and conjectured an interesting relationship between $\mathcal{S}(\Sigma_{g,k};\mathbb{Q}(q^{\frac{1}{2}}))$ and the mapping class group of $\Sigma_{g,k}$.
As related works, the structure of the skein algebra at roots of unity were investigated by Frohman and Kania-Bartoszynska \cite{FK18} and Frohman et. al. \cite{FKL21}; a presentation for the {\it stated skein algebra} was deduced by Korinman \cite{Ko23}.

We aim at uncovering the structure of $\mathcal{S}_n:=\mathcal{S}(\Sigma_{0,n+1};R)$, for any $R$, any $n\ge 4$. When seeking a presentation, it is relatively easy to find a set of relations, but hard to show that they generate the ideal $\mathcal{I}$ of defining relations. We manage to give an upper bound on the degrees of relations to generate $\mathcal{I}$, so in principle we can determine $\mathcal{I}$ in finitely many steps.

\begin{figure}[h]
  \centering
  \includegraphics[width=5.5cm]{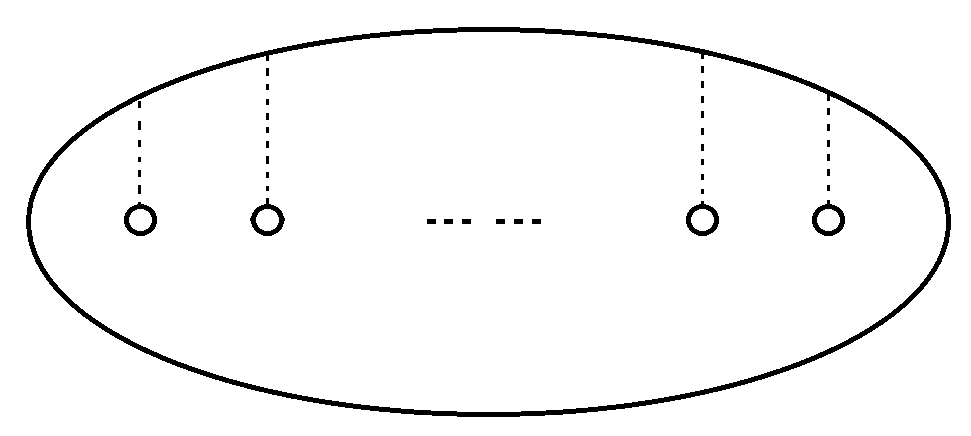}\\
  \caption{The surface $\Sigma=\Sigma_{0,n+1}$. The dotted lines are $\mathbf{z}_1,\ldots,\mathbf{z}_n$.}\label{fig:Sigma}
\end{figure}

Display $\Sigma:=\Sigma_{0,n+1}$ as a sufficiently large disk in $\mathbb{R}^2$ containing $\mathsf{p}_j:=(j,0)$, $1\le j\le n$, with a small neighborhood of $\mathsf{p}_j$ removed for each $j$ (so $\Sigma$ is compact with boundary); see Figure \ref{fig:Sigma}.
Let $\mathbf{z}_j=\{(j,y)\in\Sigma\colon y>0\}$.

\begin{figure}[h]
  \centering
  \includegraphics[width=12cm]{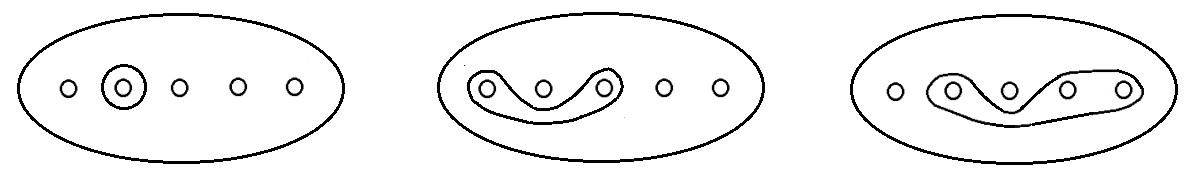}\\
  \caption{Left: $t_2$; middle: $t_{13}$; right: $t_{245}$.}\label{fig:generator}
\end{figure}

For $1\le i_1<\cdots<i_k\le n$, let $t_{i_1\cdots i_r}\in\mathcal{S}_n$ denote the element represented by a simple curve encircling $\mathsf{p}_j$ once exactly for $j\in\{i_1,\ldots,i_k\}$; see Figure \ref{fig:generator} for illustrations when $n=5$.
Let $\Delta_{i_1,\ldots,i_k}$ be the surface obtained by cutting $\Sigma$ along $\mathbf{z}_j$ for all
$j\in\{1,\ldots,n\}\setminus\{i_1,\ldots,i_k$\}. It is diffeomorphic to $\Sigma_{0,k+1}.$

Our main result is
\begin{thm}\label{thm:main}
Suppose $n\ge 4$.
\begin{enumerate}
  \item[\rm(a)] If $q+q^{-1}$ is invertible in $R$, then $\mathcal{S}_n$ is generated by
        $$\{t_{i_1\cdots i_k}\colon 1\le i_1<\cdots<i_k\le n,\ 1\le k\le 3\},$$
        and the ideal of defining relations is generated by relations of degree $\le6$ supported by $\Delta_{s_1,\ldots,s_k}$ for
        $1\le s_1<\cdots<s_k\le n$ with $k\le 6$.
  \item[\rm(b)] If $q+q^{-1}$ is not invertible in $R$, then $\mathcal{S}_n$ is generated by
        $$\{t_{i_1\cdots i_k}\colon 1\le i_1<\cdots<i_k\le n\},$$
        and the ideal of defining relations is generated by relations of degree $\le2k+2$ supported by $\Delta_{s_1,\ldots,s_k}$ for
        $1\le s_1<\cdots<s_k\le n$.
\end{enumerate}
\end{thm}

The first assertion in (a) is the genus $0$ case of \cite[Theorem 8.1]{PS00},
and the first assertion in (b) is the genus $0$ case of \cite[Theorem 1]{Bu99}.
We will reprove them.

By the phrase ``$f=0$ is a relation {\it supported} by $\Delta_{s_1,\ldots,s_k}$", we mean that $f=0$ holds in $\mathcal{S}(\Delta_{s_1,\ldots,s_k};R)\cong\mathcal{S}_k$ (so that $f=0$ holds in $\mathcal{S}_n$, due to that the inclusion $\Delta_{s_1,\ldots,s_k}\subseteq\Sigma$ induces a morphism between skein algebras).

The original problem is reduced to determining relations with degrees bounded above. Under the assumption that $q+q^{-1}$ is invertible, the problem is reduced to determining the relations of degree $\le 6$ in $\mathcal{S}_k$ for $k\le 6$. Based on Theorem \ref{thm:main}, in \cite{Ch24-1} we deduce an explicit presentation for $\mathcal{S}_n$ for each $n\ge 4$, so as to settle Problem \ref{prob:main} in genus zero case.
Setting $R=\mathbb{C}$ and $q^{\frac{1}{2}}=-1$, we recover the classical result on the structure of the trace algebra of $2\times 2$ unimodular matrices.

During the recent years, skein algebra has been found to be closely related to various structures, including quantum Teichm\"uller space \cite{BW11,Le19}, cluster algebra \cite{Mu16}, double affine Hecke algebra \cite{CS21}, higher Askey-Wilson algebra \cite{CL25}, quantized Coulomb branch \cite{AKS25,AS24}, and so on.
Our results are expected to promote understanding the various related structures, and will be beneficial to settling \cite[Question 7.1]{KMW25} which asks for presentations of some generalized skein algebras.

The content is organized as follows. In Section 2, we introduce necessary notations and conventions. In Section 3, we develop a technology of ``chopping up arcs" for simplifying links. In Section 4 we investigate the properties of the chopping-up map. In Section 5, we prove Theorem \ref{thm:main}. To improve the readability, we put the proofs of Lemma \ref{lem:reduce-multi-curve}, Lemma \ref{lem:substitution} and Lemma \ref{lem:shorten} in Section 6.

\section{Set up}

For a finite set $Y$, let $\#Y$ denote its cardinality.

Let ${\rm Sym}(r)$ denote the permutation group on $\{1,\ldots,r\}$.

Denote $q^{-1}$ as $\overline{q}$, denote $q^{-\frac{1}{2}}$ as $\overline{q}^{\frac{1}{2}}$, and so forth. Let $\alpha=q+\overline{q}$.

Use $\mathbf{a},\mathbf{b},\mathbf{x}$, etc. to denote $1$-manifolds. Let $\mathring{\mathbf{x}}$ denote the interior of $\mathbf{x}$.

Let $\pi:\Sigma\times[0,1]\to\Sigma$ denote the projection.

Let $\mathbf{z}=\cup_{k=1}^n\mathbf{z}_k$. Let $Z=\cup_{k=1}^nZ_k$, with $Z_k=\mathbf{z}_k\times[0,1]$.

When writing a sum such as $\sum_i\lambda_i\mathbf{u}_i$, we always mean a finite sum.

For $M=\Sigma\times[0,1]$ or $M=\Sigma$, an isotopy of $M$ means a continuous family $\{\varphi_t\}_{t\in[0,1]}$ of self-diffeomorphisms of $M$
such that $\varphi_0={\rm id}$ and $\varphi_t$ fixes $\partial M$ pointwise for all $t$.

\begin{defn}
\rm A self-diffeomorphism $\varphi$ of $\Sigma\times[0,1]$ fixing $\partial(\Sigma\times[0,1])$ pointwise is called a {\it congruence} if $\varphi(Z_k)=Z_k$ for all $k$ and $\pi(\varphi(\mathsf{a},z))=\pi(\varphi(\mathsf{a},0))$ for all $(\mathsf{a},z)\in\Sigma\times[0,1]$.

Given $1$-submanifolds $\mathbf{x},\mathbf{x}'\subset\Sigma\times[0,1]$, say that $\mathbf{x}$ is {\it congruent to} $\mathbf{x}'$ and denote $\mathbf{x}\cong \mathbf{x}'$, if there exists a congruence $\varphi$ such that $\varphi(\mathbf{x})=\mathbf{x}'$.
\end{defn}

As usual, we present links in $\Sigma\times(0,1)$ via projection diagrams in $\Sigma$.
Given a link $\mathbf{l}\subset\Sigma\times(0,1)$, let $[\mathbf{l}]\in\mathcal{S}_n$ denote the element represented by $\mathbf{l}$.
Draw $\mathbf{l}_1$ above $\mathbf{l}_2$ when defining $[\mathbf{l}_1][\mathbf{l}_2]$ in $\mathcal{S}_n$.

Unless otherwise specified, a $1$-submanifold $\mathbf{x}\subset\Sigma\times[0,1]$ is always assumed to be compact and {\it in generic position}, in the sense that up to diffeomorphism, $\pi(\mathbf{x})$ is stable under small perturbations. In particular, the following holds: $\mathbf{x}$ intersects $Z$ transversally; $\partial \mathbf{x}\cap Z=\emptyset$; $\#\pi^{-1}(\mathsf{a})=1$ for all $\mathsf{a}\in\pi(\mathbf{x})\setminus{\rm Cr}(\mathbf{x})$, where ${\rm Cr}(\mathbf{x})$ is a finite subset of $\pi(\mathring{\mathbf{x}})\setminus\mathbf{z}$ such that $\#\pi^{-1}(\mathsf{c})=2$ for each $\mathsf{c}\in{\rm Cr}(\mathbf{x})$.
Furthermore, we always assume that each connected component of $\mathbf{x}$ diffeomorphic to $S^1$ lies in $\Sigma\times(0,1)$.

Each $\mathsf{c}\in{\rm Cr}(\mathbf{x})$ is caused by a crossing made by two short arcs $\mathbf{a}_{\pm}\subset\mathbf{x}$;
let ${\rm over}(\mathsf{c})\in\mathbf{a}_+$, ${\rm under}(\mathsf{c})\in\mathbf{a}_-$ respectively denote the upper and lower point that constitute $\pi^{-1}(\mathsf{c})$. Abusing the notation, we also call $\mathsf{c}$ a crossing.

Let ${\rm cn}(\mathbf{x})=\#{\rm Cr}(\mathbf{x})$. Call $\mathbf{x}$ {\it simple} if ${\rm cn}(\mathbf{x})=0$; in this case, we may identify $\mathbf{x}$ with $\pi(\mathbf{x})$ which is a $1$-submanifold of $\Sigma$.

Let $|\mathbf{x}|_i=\#(\mathbf{x}\cap Z_i)$ for $1\le i\le n$. Let ${\rm supp}(\mathbf{x})=\{i\colon |\mathbf{x}|_i>0\}$.
Define the {\it degree} of $\mathbf{x}$ as $|\mathbf{x}|:={\sum}_{i=1}^n|\mathbf{x}|_i=\#(\mathbf{x}\cap Z).$
Let $\|\mathbf{x}\|=(|\mathbf{x}|,{\rm cn}(\mathbf{x}))$.

Introduce a linear order $\prec$ on $\mathbb{Z}_{\ge0}^2$, by declaring $(m',c')\prec(m,c)$ for $m'<m$, and $(m,c')\prec(m,c)$ for $c'<c$. Denote $(m',c')\preceq(m,c)$ if $(m',c')\prec(m,c)$ or $(m',c')=(m,c)$.
Given $\mathbf{x},\mathbf{x}'\subset\Sigma\times[0,1]$, when $\|\mathbf{x}\|\prec\|\mathbf{x}'\|$, we say that $\mathbf{x}$ is {\it simpler} than $\mathbf{x}'$, and abusively denote $\mathbf{x}\prec\mathbf{x}'$.

\begin{defn}
\rm A {\it stacked link} is a disjoint union $\mathbf{u}=\mathbf{k}_1\sqcup\cdots\sqcup \mathbf{k}_r$ such that $\mathbf{k}_i$ is a knot in $\Sigma\times(z_i,z_{i-1})$, where $1=z_0>\cdots>z_r=0$.

Define an equivalence relation $\sim$ among stacked links with equally many components, by declaring
$\mathbf{u}\sim\mathbf{u}'$ if $({\rm id}_\Sigma\times\rho)(\mathbf{u})=\mathbf{u}'$ for some orientation-preserving diffeomorphism $\rho:[0,1]\to[0,1]$. Note that if $\mathbf{u}\sim\mathbf{u}'$, then $\|\mathbf{u}\|=\|\mathbf{u}'\|$.
Ignoring the information of $z_i$'s, we denote $\mathbf{k}_1\sqcup\cdots\sqcup \mathbf{k}_r$ as $\mathbf{k}_1\cdots \mathbf{k}_r$, which is well-defined up to $\sim$.
\end{defn}

Let $\mathcal{F}$ denote the quotient of the free $R$-module generated by equivalence classes of stacked links by the submodule generated by elements of the forms $\mathbf{u}\mathbf{o}+\alpha\mathbf{u}$ and $\mathbf{o}\mathbf{u}+\alpha\mathbf{u}$, where $\mathbf{u}$ is a stacked link and $\mathbf{o}$ is a simple curve of degree $0$.
It is a $R$-algebra with multiplication defined via stacking.
Let $\tilde{\theta}:\mathcal{F}\to\mathcal{S}_n$ denote the $R$-algebra morphism sending a stacked link $\mathbf{u}$ to $[\mathbf{u}]$.
For $a_1,a_2\in\mathcal{F}$, we say ``$a_1=a_2$ in $\mathcal{S}_n$" if $\tilde{\theta}(a_1)=\tilde{\theta}(a_2)$.

Given $1$-submanifolds $\mathbf{x},\mathbf{y}$ with $\mathring{\mathbf{x}}\cap\mathring{\mathbf{y}}=\emptyset$, let
$${\rm Cr}(\mathbf{x},\mathbf{y})=\big\{\mathsf{c}\in{\rm Cr}(\mathbf{x}\cup \mathbf{y})\colon{\rm over}(\mathsf{c})\in \mathbf{x},\
{\rm under}(\mathsf{c})\in \mathbf{y}\big\}.$$

Call a $1$-submanifold $\mathbf{c}$ an {\it arc} if it is diffeomorphic to $[0,1]$.

Given an oriented arc $\mathbf{c}$, write $\partial \mathbf{c}=\{\partial_-\mathbf{c},\partial_+\mathbf{c}\}$, so that $\mathbf{c}$ is oriented from $\partial_-\mathbf{c}$ to $\partial_+\mathbf{c}$. Starting at $\partial_-\mathbf{c}$, walk along $\mathbf{c}$ towards $\partial_+\mathbf{c}$,
record $x_i$ (resp. $x_i^{-1}$) whenever passing through $\mathbf{z}_i$ from left to right (resp. from right to left).
If all the recordings are $x_{i_1}^{\nu_1},\ldots,x_{i_m}^{\nu_m}$, then we put
${\rm word}(\mathbf{c})=x_{i_1}^{\nu_1}\cdots x_{i_m}^{\nu_m}$.

Let $\mathbf{a},\mathbf{b}$ be arcs with $\mathring{\mathbf{a}}\cap\mathring{\mathbf{b}}=\emptyset$.
When $\#(\partial\mathbf{a}\cap\partial\mathbf{b})=1$, let
$\mathbf{a}\mathbf{b}=\mathbf{a}\cup \mathbf{b}$; when $\partial\mathbf{a}=\partial\mathbf{b}$, let ${\rm tr}(\mathbf{a}\mathbf{b})=\mathbf{a}\cup \mathbf{b}$ which is a knot. We may perturb $\mathbf{a}$ or $\mathbf{b}$ if necessary, to make $\mathbf{a}\cup \mathbf{b}$ in generic position. Moreover, we require that the tangent vector of $\mathbf{a}$ is parallel to that of $\mathbf{b}$ at each joining point. Such conventions will be always adopted.

Let $\mathbf{x}$ be a $1$-submanifold. When $\mathbf{a}\subset \mathbf{x}$ is an arc with $\pi(\partial \mathbf{a})\cap{\rm Cr}(\mathbf{x})=\emptyset$, we call it an arc of $\mathbf{x}$, and denote $\langle\mathbf{x}|\mathbf{a}\rangle$ for $\mathbf{x}\setminus\mathring{\mathbf{a}}$. Given another arc $\mathbf{c}$ with $\mathbf{c}\cap\langle\mathbf{x}|\mathbf{a}\rangle=\partial\mathbf{c}=\partial\mathbf{a}$, let
$(\mathbf{x}|\mathbf{a}|\mathbf{c})=\langle\mathbf{x}|\mathbf{a}\rangle\cup\mathbf{c}$.

Let ${\rm Ar}_d(\mathbf{x})$ denote the set of degree $d$ arcs of $\mathbf{x}$. Let ${\rm Ar}(\mathbf{x})=\cup_{d\ge 0}{\rm Ar}_d(\mathbf{x})$.

Given $\mathsf{a},\mathsf{b}\in\Sigma\times\{z\}$ with $z\in\{0,1\}$, let $F_z(\mathsf{a},\mathsf{b})$ denote the set of $1$-submanifolds $\mathbf{x}=\mathbf{c}\sqcup\mathbf{l}$, where $\mathbf{l}$ is a link and $\mathbf{c}$ is an arc with $\partial\mathbf{c}=\{\mathsf{a},\mathsf{b}\}$.
Let $\mathcal{S}(\mathsf{a},\mathsf{b})$ denote the $R$-module generated by isotopy classes of elements of $F_z(\mathsf{a},\mathsf{b})$, modulo skein relations similar to (\ref{eq:local}).
Let $[\mathbf{x}]\in\mathcal{S}(\mathsf{a},\mathsf{b})$ denote the element represented by $\mathbf{x}$.
When $z=0$ (resp. $z=1$), $\mathcal{S}(\mathsf{a},\mathsf{b})$ is a left (resp. right) $\mathcal{S}_n$-module.
This is a special case of the {\it relative skein module} proposed in \cite[Definition 3.8]{Pr99} and \cite[Definition 5.1]{PS00}.

\begin{rmk}\label{rmk:explain}
\rm We pause to give explanations for some notions.

(i) The notion ``congruence" is introduced for keeping track of the information about $\mathbf{x}\cap Z$ and the crossings of $\mathbf{x}$,
for a 1-manifold $\mathbf{x}$. Clearly, $\|\mathbf{x}\|=\|\mathbf{x}'\|$ if $\mathbf{x}\cong \mathbf{x}'$.

(ii) 
Suppose $\mathbf{s}_1\sqcup\cdots\sqcup \mathbf{s}_r\subset\Sigma$ is a multicurve, where each $\mathbf{s}_i$ is a nonnullhomotopic simple curve.
Given $\sigma\in{\rm Sym}(r)$, we can construct a stacked link $\mathbf{s}_{\sigma(1)}\cdots \mathbf{s}_{\sigma(r)}$. Be aware that $\mathbf{s}_{\sigma(1)}\cdots \mathbf{s}_{\sigma(r)}\nsim \mathbf{s}_{\tau(1)}\cdots \mathbf{s}_{\tau(r)}$ if $\sigma\ne\tau$.

(iii) When $\mathbf{a}\in{\rm Ar}(\mathbf{x})$ with ${\rm Cr}(\langle\mathbf{x}|\mathbf{a}\rangle,\mathbf{a})=\emptyset$, we can take a congruence $\varphi$ such that if $\mathbf{a}':=\varphi(\mathbf{a})$, $\mathbf{x}'=\varphi(\mathbf{x})$, then
$$\mathbf{a}'\subset\Sigma\times[1/2,1], \qquad \langle\mathbf{x}'|\mathbf{a}'\rangle\subset\Sigma\times[0,1/2], \qquad  \mathbf{x}'\cap\Sigma\times\{1/2\}=\partial \mathbf{a}'.$$
Let $\{\mathsf{a}_{\pm}\}\subset\Sigma\times\{0\}$ be the image of $\partial\mathbf{a}'$ under the map $\Sigma\times[\frac{1}{2},1]\to\Sigma\times[0,1]$, $(\mathsf{x},v)\mapsto(\mathsf{x},2v-1)$. Then $\mathbf{a}'$ corresponds to an arc $\mathbf{a}''\in F_0(\mathsf{a}_-,\mathsf{a}_+)$.
On this account, we may bear in mind that any operation applied to $\mathbf{a}''$ is applicable to $\mathbf{a}'$, and also to $\mathbf{a}$.

The situation is similar when ${\rm Cr}(\mathbf{a},\langle\mathbf{x}|\mathbf{a}\rangle)=\emptyset$.
\end{rmk}

\begin{defn}
\rm Let $\mathbf{x}$ be a $1$-submanifold of $\Sigma$. Call $\mathbf{a}\in {\rm Ar}_2(\mathbf{x})$ {\it shrinkable} if $\mathbf{a}$ can be {\it reduced}, in the sense that it can be isotoped into some degree $0$ arc $\mathbf{b}$ through an isotopy fixing $\langle\mathbf{x}|\mathbf{a}\rangle$, so that $(\mathbf{x}|\mathbf{a}|\mathbf{b})$ is isotopic to $\mathbf{x}$. Call $\mathbf{x}$ {\it reducible} if it has a shrinkable arc; otherwise call $\mathbf{x}$ {\it irreducible}.
\end{defn}

Observe that if $\mathbf{x}$ is reducible, then its shrinkable arcs can be reduced successively.
Also observe that $\mathbf{x}$ is reducible if and only if ${\rm word}(\mathbf{c})$ is reducible for some oriented arc $\mathbf{c}\subseteq\mathbf{x}$.

\begin{defn}
\rm Call isotopy $\phi_t:\Sigma\to\Sigma$ {\it fine} if $\phi_t(\mathbf{z}_k)=\mathbf{z}_k$ for all $k,t$.
\end{defn}

\begin{lem}\label{lem:reduce-multi-curve}
If two irreducible multicurves $\mathbf{m}_0,\mathbf{m}_1$ are isotopic, then there exists a fine isotopy $\varphi_t$ with $\varphi_1(\mathbf{m}_0)=\mathbf{m}_1$.
\end{lem}

Let $\mathcal{V}$ denote the free $R$-module generated by fine isotopy classes of irreducible multicurves, which, by Lemma \ref{lem:reduce-multi-curve}, can be identified with the free $R$-module generated by isotopy classes of multicurves.

For a link $\mathbf{l}$, let $\Theta'(\mathbf{l})$ denote the $R$-linear combination of fine isotopy classes of (possibly reducible) multicurves obtained by resolving all crossings of $\mathbf{l}$, and let $\Theta(\mathbf{l})\in\mathcal{V}$ denote the $R$-linear combination of fine isotopy classes of irreducible multicurves obtained from $\Theta'(\mathbf{l})$ by reducing shrinkable arcs and removing degree $0$ simple curves (and multiplying by a power of $-\alpha$), whenever possible.
We emphasize that $\Theta'(\mathbf{l})$, $\Theta(\mathbf{l})$ are independent of the order of resolving crossings.

\begin{nota}\label{nota:skein-triple}
\rm Given a link $\mathbf{l}$ and $\mathsf{c}\in{\rm Cr}(\mathbf{l})$, let $\mathbf{l}^{\mathsf{c}}_\infty,\mathbf{l}^{\mathsf{c}}_0$ denote the links obtained by resolving $\mathsf{c}$ such that $(\mathbf{l},\mathbf{l}^{\mathsf{c}}_\infty,\mathbf{l}^{\mathsf{c}}_0)$ is a skein triple (as in Figure \ref{fig:local}).
\end{nota}
By definition, each skein triple has the form $(\mathbf{l},\mathbf{l}^{\mathsf{c}}_\infty,\mathbf{l}^{\mathsf{c}}_0)$, for some link $\mathbf{l}$ and some $\mathsf{c}\in{\rm Cr}(\mathbf{l})$, and $\mathbf{l}=q^{\frac{1}{2}}\mathbf{l}^{\mathsf{c}}_\infty+\overline{q}^{\frac{1}{2}}\mathbf{l}^{\mathsf{c}}_0$ in $\mathcal{S}_n$.

For each skein triple $(\mathbf{l},\mathbf{l}^{\mathsf{c}}_\infty,\mathbf{l}^{\mathsf{c}}_0)$,
we have
$\Theta'(\mathbf{l})=q^{\frac{1}{2}}\Theta'(\mathbf{l}^{\mathsf{c}}_\infty)+\overline{q}^{\frac{1}{2}}\Theta'(\mathbf{l}^{\mathsf{c}}_0).$
To see this, just resolve $\mathsf{c}$ first for $\mathbf{l}$, then the remaining crossings of $\mathbf{l}$ coincide with the crossings of $\mathbf{l}^{\mathsf{c}}_\infty$ and $\mathbf{l}^{\mathsf{c}}_0$, so we can resolve them synchronously.
Hence
$$\Theta(\mathbf{l})=q^{\frac{1}{2}}\Theta(\mathbf{l}^{\mathsf{c}}_\infty)+\overline{q}^{\frac{1}{2}}\Theta(\mathbf{l}^{\mathsf{c}}_0).$$

Consequently, $\Theta$ induces a map $\mathcal{S}_n\to\mathcal{V}$, which is abusively denoted by $\Theta$, too.
By \cite[Corollary 4.1]{SW07}, it is an isomorphism of $R$-modules.

\section{Chopping up arcs}

Let $\mathcal{T}$ denote the $R$-subalgebra of $\mathcal{F}$ generated by
$$G:=\begin{cases} \{t_{i_1\cdots i_k}\colon 1\le i_1<\cdots<i_k\le n,\ 1\le k\le 3\},&\alpha^{-1}\in R \\
\{t_{i_1\cdots i_k}\colon 1\le i_1<\cdots<i_k\le n\},&\alpha^{-1}\notin R\end{cases}.$$
Let $\theta:\mathcal{T}\to\mathcal{S}_n$ denote the restriction of $\tilde{\theta}$ to $\mathcal{T}$.

\begin{defn}
\rm Call an arc $\mathbf{a}$ {\it unshortenable} if $|\mathbf{a}|=\#{\rm supp}(\mathbf{a})$, and {\it shortenable} if $|\mathbf{a}|>\#{\rm supp}(\mathbf{a})$. Call $\mathbf{a}$ {\it minimal shortenable} if $|\mathbf{a}|=\#{\rm supp}(\mathbf{a})+1$.
\end{defn}

If $|\mathbf{a}|\ge n+1$, then $\mathbf{a}$ must be shortenable.

With an orientation chosen, suppose
${\rm word}(\mathbf{a})=x_{i_1}^{\nu_1}\cdots x_{i_m}^{\nu_m}$, then $\mathbf{a}$ is minimal shortenable if and only if $i_1,\ldots,i_{m-1}$ are distinct and $i_m=i_1$.

\begin{nota}
\rm When $\alpha^{-1}\in R$, let $h=3$, and let ${\rm Ar}^\ast(\mathbf{x})={\rm Ar}_3(\mathbf{x})$, for a $1$-submanifold $\mathbf{x}\subset\Sigma\times[0,1]$.
When $\alpha^{-1}\not\in R$, let $h=n+1$, and let ${\rm Ar}^\ast(\mathbf{x})$ denote the set of minimal shortenable arcs of $\mathbf{x}$.
\end{nota}

\begin{lem}\label{lem:substitution}
Suppose $\alpha^{-1}\in R$. Given $\mathsf{a},\mathsf{b}\in\Sigma\times\{z\}$ with $z\in\{0,1\}$, let $P(\mathsf{a},\mathsf{b})$ denote the subset of $F_z(\mathsf{a},\mathsf{b})$ consisting of simple irreducible arcs $\mathbf{c}$ with $|\mathbf{c}|\le 2$.
Let $\mathbf{a}\in F_z(\mathsf{a},\mathsf{b})$ be an arbitrary degree $3$ arc.
\begin{enumerate}
  \item[\rm(i)] If $z=0$, then $[\mathbf{a}]={\rm ch}_u(\mathbf{a}):=\sum_sa_s[\mathbf{c}_s]$ in $\mathcal{S}(\mathsf{a},\mathsf{b})$ for some
        $a_s\in\mathcal{T}$, $\mathbf{c}_s\in P(\mathsf{a},\mathsf{b})$ such that $|a_s|_i+|\mathbf{c}_s|_i\le |\mathbf{a}|_i$ for all $s,i$.
  \item[\rm(ii)] If $z=1$, then $[\mathbf{a}]={\rm ch}_d(\mathbf{a}):=\sum_t[\mathbf{d}_t]b_t$ in $\mathcal{S}(\mathsf{a},\mathsf{b})$ for some
        $b_t\in\mathcal{T}$, $\mathbf{d}_t\in P(\mathsf{a},\mathsf{b})$ such that $|b_t|_i+|\mathbf{d}_t|_i\le |\mathbf{a}|_i$ for all $t,i$.
\end{enumerate}
\end{lem}

\begin{lem}\label{lem:shorten}
Suppose $\alpha^{-1}\not\in R$. Given $\mathsf{a},\mathsf{b}\in\Sigma\times\{z\}$ with $z\in\{0,1\}$, let $Q(\mathsf{a},\mathsf{b})$ denote the subset of $F_z(\mathsf{a},\mathsf{b})$ consisting of simple irreducible unshortenable arcs.
Let $\mathbf{a}\in F_z(\mathsf{a},\mathsf{b})$ be an arbitrary minimal shortenable arc.
\begin{enumerate}
  \item[\rm(i)] If $z=0$, then $[\mathbf{a}]={\rm ch}_u(\mathbf{a}):=\sum_sa_s[\mathbf{c}_s]$ in $\mathcal{S}(\mathsf{a},\mathsf{b})$ for some
        $a_s\in\mathcal{T}$ and $\mathbf{c}_s\in Q(\mathsf{a},\mathsf{b})$ such that $|a_s|_i+|\mathbf{c}_s|_i\le|\mathbf{a}|_i$ for all $s,i$.
  \item[\rm(ii)] If $z=1$, then $[\mathbf{a}]={\rm ch}_d(\mathbf{a}):=\sum_v[\mathbf{d}_t]b_t$ in $\mathcal{S}(\mathsf{a},\mathsf{b})$ for some
        $b_t\in\mathcal{T}$ and $\mathbf{d}_t\in Q(\mathsf{a},\mathsf{b})$ such that $|b_t|_i+|\mathbf{d}_t|_i\le|\mathbf{a}|_i$ for all $t,i$.
\end{enumerate}
\end{lem}

\begin{nota}
\rm Given a knot $\mathbf{k}$ and $\mathsf{c}\in{\rm Cr}(\mathbf{k})$, let $\mathbf{k}^\star(\mathsf{c})$ denote the knot obtained by changing the type of $\mathsf{c}$.
Recalling Notation \ref{nota:skein-triple}, note that one of $\mathbf{k}^{\mathsf{c}}_\infty,\mathbf{k}^{\mathsf{c}}_0$ is a two-component link, and the other is a knot; let $\mathbf{k}^\flat(\mathsf{c})$ denote the knot.
Put
$$\epsilon(\mathsf{c})=
\begin{cases}
1,&\mathbf{k}^\flat(\mathsf{c})=\mathbf{k}^{\mathsf{c}}_0 \\
-1,&\mathbf{k}^\flat(\mathsf{c})=\mathbf{k}^{\mathsf{c}}_\infty
\end{cases};  \qquad
\delta(\mathsf{c})=\overline{q}^{\frac{1}{2}\epsilon(\mathsf{c})}-q^{\frac{3}{2}\epsilon(\mathsf{c})}.$$
\end{nota}

\begin{figure}[h]
  \centering
  \includegraphics[width=11.8cm]{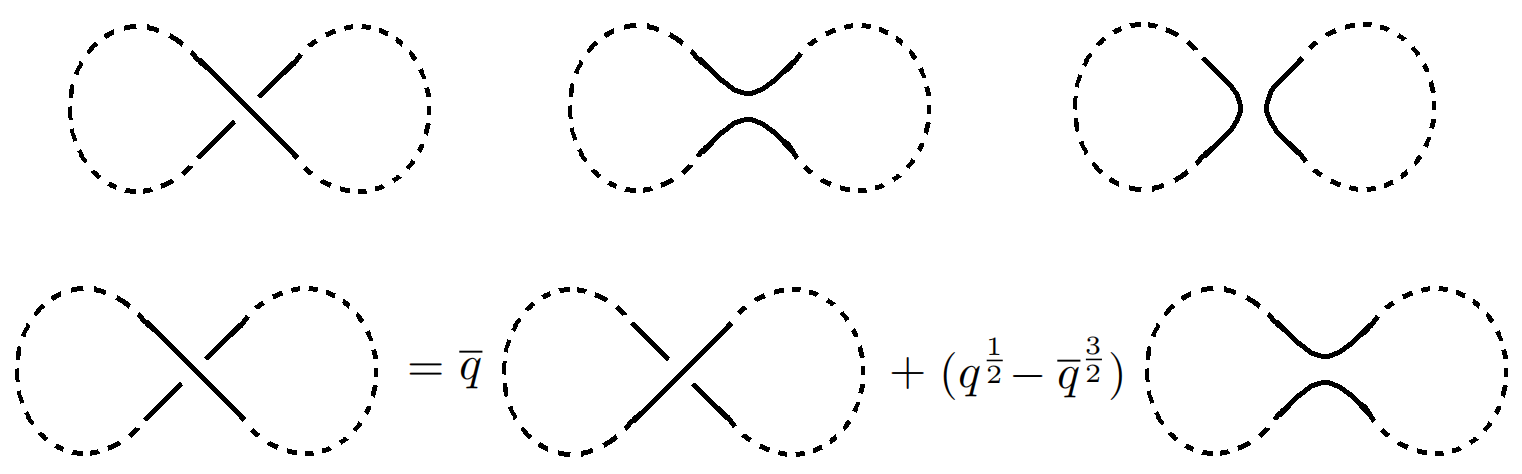}\\
  \caption{First row (from left to right): a crossing $\mathsf{c}$ of a knot $\mathbf{k}$; $\mathbf{k}^{\mathsf{c}}_\infty$;
  $\mathbf{k}^{\mathsf{c}}_0$. Second row: $\mathbf{k}=q^{\epsilon(\mathsf{c})}\mathbf{k}^\star(\mathsf{c})+\delta(\mathsf{c})\mathbf{k}^\flat(\mathsf{c})$, which holds in $\mathcal{S}_n$;
  here $\epsilon(\mathsf{c})=-1$, $\delta(\mathsf{c})=q^{\frac{1}{2}}-\overline{q}^{\frac{3}{2}}$.
  The dotted arcs stand for the remaining parts, which are possibly very complicated but irrelevant.}\label{fig:change}
\end{figure}

In $\mathcal{S}_n$, since
$\mathbf{k}=q^{\frac{1}{2}}\mathbf{k}^{\mathsf{c}}_\infty+\overline{q}^{\frac{1}{2}}\mathbf{k}^{\mathsf{c}}_0$ and $\mathbf{k}^\star(\mathsf{c})=\overline{q}^{\frac{1}{2}}\mathbf{k}^{\mathsf{c}}_\infty+q^{\frac{1}{2}}\mathbf{k}^{\mathsf{c}}_0$,
we have
\begin{align}
\mathbf{k}=q^{\epsilon(\mathsf{c})}\mathbf{k}^\star(\mathsf{c})+\delta(\mathsf{c})\mathbf{k}^\flat(\mathsf{c}).  \label{eq:K}
\end{align}
See Figure \ref{fig:change} for illustration.

Let $\mathbf{a}\in{\rm Ar}(\mathbf{k})$, with an orientation $\o$ chosen. Suppose
\begin{align*}
{\rm Cr}(\langle \mathbf{k}|\mathbf{a}\rangle,\mathbf{a})=\{\mathsf{c}^1,\ldots,\mathsf{c}^k\},  \qquad
{\rm Cr}(\mathbf{a},\langle \mathbf{k}|\mathbf{a}\rangle)=\{\mathsf{c}_1,\ldots,\mathsf{c}_\ell\},
\end{align*}
both listed in the order determined by $\o$.
Be aware that the crossings from ${\rm Cr}(\mathbf{a})$ are not included.
Put $\mathbf{k}^{[0]}=\mathbf{k}_{[0]}=\mathbf{k}$; recursively put
\begin{align*}
\mathbf{k}^{[i]}&=(\mathbf{k}^{[i-1]})^\star(\mathsf{c}^i),  \qquad 1\le i\le k;  \\
\mathbf{k}_{[i]}&=(\mathbf{k}_{[i-1]})^\star(\mathsf{c}_i),  \qquad 1\le i\le \ell.
\end{align*}
In other words, $\mathbf{k}^{[i]}$ (resp. $\mathbf{k}_{[i]}$) is the knot obtained from $\mathbf{k}$ by changing the types of $\mathsf{c}^1,\ldots,\mathsf{c}^i$ (resp. $\mathsf{c}_1,\ldots,\mathsf{c}_i$).
Let $\mathbf{k}^{\mathbf{a}}=\mathbf{k}^{[k]}$, $\mathbf{k}_{\mathbf{a}}=\mathbf{k}_{[\ell]}$.
Intuitively, $\mathbf{k}^{\mathbf{a}}$ (resp. $\mathbf{k}_{\mathbf{a}}$) is obtained from $\mathbf{k}$ by pulling $\mathbf{a}$ up to the top (resp. pushing $\mathbf{a}$ down to the bottom) via crossing-changes. Clearly, $\|\mathbf{k}^{\mathbf{a}}\|=\|\mathbf{k}_{\mathbf{a}}\|=\|\mathbf{k}\|$.
Let $\mathbf{a}^\diamond\in{\rm Ar}(\mathbf{k}^{\mathbf{a}})$, $\mathbf{a}_\diamond\in{\rm Ar}(\mathbf{k}_{\mathbf{a}})$ denote the arcs inherited from $\mathbf{a}$.

By (\ref{eq:K}), in $\mathcal{S}_n$,
\begin{align*}
\mathbf{k}^{[i-1]}&=q^{\epsilon(\mathsf{c}^i)}\mathbf{k}^{[i]}+\delta(\mathsf{c}^i)(\mathbf{k}^{[i-1]})^\flat(\mathsf{c}^i),  \\
\mathbf{k}_{[i-1]}&=q^{\epsilon(\mathsf{c}_i)}\mathbf{k}_{[i]}+\delta(\mathsf{c}_i)(\mathbf{k}_{[i-1]})^\flat(\mathsf{c}_i).
\end{align*}
Hence
\begin{align}
\mathbf{k}=q^{\hat{\epsilon}(\mathbf{k},\mathbf{a})}\mathbf{k}^{\mathbf{a}}+\mathfrak{r}_u(\mathbf{k},\mathbf{a})
=q^{\check{\epsilon}(\mathbf{k},\mathbf{a})}\mathbf{k}_{\mathbf{a}}+\mathfrak{r}_d(\mathbf{k},\mathbf{a})  \qquad \text{in}\ \mathcal{S}_n, \label{eq:equality-in-Sn}
\end{align}
where $\hat{\epsilon}(\mathbf{k},\mathbf{a})={\sum}_{i=1}^k\epsilon(\mathsf{c}^i)$, $\check{\epsilon}(\mathbf{k},\mathbf{a})={\sum}_{i=1}^\ell\epsilon(\mathsf{c}_i)$, and
\begin{align*}
\mathfrak{r}_u(\mathbf{k},\mathbf{a})&={\sum}_{i=1}^k
q^{\sum_{t=1}^{i-1}\epsilon(\mathsf{c}^t)}\delta(\mathsf{c}^i)\cdot(\mathbf{k}^{[i-1]})^\flat(\mathsf{c}^i), \\ 
\mathfrak{r}_d(\mathbf{k},\mathbf{a})&={\sum}_{i=1}^\ell q^{\sum_{t=1}^{i-1}\epsilon(\mathsf{c}_t)}\delta(\mathsf{c}_i)\cdot(\mathbf{k}_{[i-1]})^\flat(\mathsf{c}_i).  
\end{align*}

\begin{rmk}
\rm The precise expression of $\mathfrak{r}_u(\mathbf{k},\mathbf{a})$ is not important; just remember that it is a $R$-linear combination of knots simpler than $\mathbf{k}$. Similarly for $\mathfrak{r}_d(\mathbf{k},\mathbf{a})$.
\end{rmk}

\begin{nota}
\rm For $m\ge 2h+1$ and $c\ge 0$, let $\mathcal{L}_{m,c}$ denote the $R$-submodule of $\mathcal{F}$ generated by elements of the following three types:
\begin{enumerate}
  \item[\rm(i)] $\sum_ia_i\mathbf{k}_i\in\tilde{\theta}^{-1}(0)$, where $a_i\in R$ and $\mathbf{k}_i$ is a knot with $\|\mathbf{k}_i\|\prec(m,c)$;
  \item[\rm(ii)] $\mathbf{k}-q^{\hat{\epsilon}(\mathbf{k},\mathbf{a})}\sum_sa_s(\mathbf{k}^{\mathbf{a}}|\mathbf{a}^\diamond|\mathbf{c}_s)
        -\mathfrak{r}_u(\mathbf{k},\mathbf{a})$ for a knot $\mathbf{k}$ with $\|\mathbf{k}\|\prec(m,c)$ and oriented
        $\mathbf{a}\in{\rm Ar}^\ast(\mathbf{k})$ such that ${\rm ch}_u(\mathbf{a}^\diamond)=\sum_sa_s[\mathbf{c}_s]$;
  \item[\rm(iii)] $\mathbf{k}-q^{\check{\epsilon}(\mathbf{k},\mathbf{a})}\sum_t(\mathbf{k}_{\mathbf{a}}|\mathbf{a}_\diamond|\mathbf{d}_t)b_t
        -\mathfrak{r}_d(\mathbf{k},\mathbf{a})$ for a knot $\mathbf{k}$ with $\|\mathbf{k}\|\prec(m,c)$ and oriented
        $\mathbf{a}\in{\rm Ar}^\ast(\mathbf{k})$ such that ${\rm ch}_d(\mathbf{a}_\diamond)=\sum_t[\mathbf{d}_t]b_t$;
  \item[\rm(iv)] $\mathbf{j}(\sum_ia_i\mathbf{l}_i)\mathbf{j}'$, where $a_i\in R$, $\mathbf{l}_i$ is a stacked link with $|\mathbf{l}_i|<m$,
        $|\mathbf{j}|+|\mathbf{l}_i|+|\mathbf{j}'|\le m$, and $\sum_ia_i\mathbf{l}_i\in\tilde{\theta}^{-1}(0)$.
\end{enumerate}
Note that for any $0\le c'<c$, type (iv) elements of $\mathcal{L}_{m,c}$ coincide with those of $\mathcal{L}_{m,c'}$.

For $m\le 2h$ and $c\ge 0$, let $\mathcal{L}_{m,c}$ denote the $R$-submodule of $\mathcal{F}$ generated by elements of the form $\sum_ia_i\mathbf{l}_i$ such that $a_i\in R$, and $\mathbf{l}_i$ is a stacked link with $|\mathbf{l}_i|\le 2h$, and $\sum_ia_i[\mathbf{l}_i]=0$.

Given a stacked link $\mathbf{l}$ with $\|\mathbf{l}\|=(m,c)$, let $\mathcal{L}_{\mathbf{l}}=\mathcal{L}_{m,c}$.
\end{nota}

\begin{nota}
\rm For $v,z\in\mathcal{F}$, use $v\equiv z\pmod{\mathcal{L}_{m,c}}$ to indicate $v-z\in\mathcal{L}_{m,c}$.
\end{nota}

Suppose $\mathbf{a}\in{\rm Ar}^\ast(\mathbf{k})$, and ${\rm ch}_u(\mathbf{a}^\diamond)=\sum_sa_s[\mathbf{c}_s]$,
${\rm ch}_d(\mathbf{a}_\diamond)=\sum_t[\mathbf{d}_t]b_t$.
Put
\begin{alignat}{2}
\varepsilon_u(\mathbf{k},\mathbf{a})&=\mathfrak{r}_u(\mathbf{k},\mathbf{a})
+q^{\hat{\epsilon}(\mathbf{k},\mathbf{a})}{\sum}_sa_s\cdot(\mathbf{k}^{\mathbf{a}}|\mathbf{a}^\diamond|\mathbf{c}_s)
&&\pmod{\mathcal{L}_{\mathbf{k}}},  \label{eq:substitution-1}  \\
\varepsilon_d(\mathbf{k},\mathbf{a})&=\mathfrak{r}_d(\mathbf{k},\mathbf{a})
+q^{\check{\epsilon}(\mathbf{k},\mathbf{a})}{\sum}_t(\mathbf{k}_{\mathbf{a}}|\mathbf{a}_\diamond|\mathbf{d}_t)\cdot b_t
&&\pmod{\mathcal{L}_{\mathbf{k}}}. \label{eq:substitution-2}
\end{alignat}

\begin{rmk}\label{rmk:convention}
\rm (i) We could have defined $\varepsilon_u(\mathbf{k},\mathbf{a}),\varepsilon_d(\mathbf{k},\mathbf{a})$ as elements of $\mathcal{F}/\mathcal{L}_{\mathbf{k}}$. However, as will be seen below, it is usually necessary to keep track of explicit representatives.
So it is more convenient to adopt notations such as $w\equiv w'\pmod{\mathcal{L}_{\mathbf{k}}}$, to indicate $w=w'$ in $\mathcal{F}/\mathcal{L}_{\mathbf{k}}$.

(ii) For each $s$, the congruence class of $(\mathbf{k}^{\mathbf{a}}|\mathbf{a}^\diamond|\mathbf{c}_s)$ depends on the choice of the concrete arc $\mathbf{c}_s$, but the isotopy class does not.
Moreover, $\mathfrak{r}_u(\mathbf{k},\mathbf{a})$ depends on the orientation of $\mathbf{a}$, but it is a $R$-linear combination of knots simpler than $\mathbf{k}$ and $\tilde{\theta}(\mathfrak{r}_u(\mathbf{k},\mathbf{a}))$ depends only on $\mathbf{k},\mathbf{a}$.
Thus, $\varepsilon_u(\mathbf{k},\mathbf{a})$ is well-defined, due to the definition of $\mathcal{L}_{\mathbf{k}}$. Similarly for $\varepsilon_d(\mathbf{k},\mathbf{a})$.

It is clear that
$$\tilde{\theta}(\varepsilon_u(\mathbf{k},\mathbf{a}))=\tilde{\theta}(\varepsilon_d(\mathbf{k},\mathbf{a}))=[\mathbf{k}].$$

(iii) Sometimes we will rewrite (\ref{eq:substitution-1}) as $\varepsilon_u(\mathbf{k},\mathbf{a})=\sum_ic_i\mathbf{p}_i$, meaning that for each $i$, either $c_i=a_s$, $\mathbf{p}_i=(\mathbf{k}^{\mathbf{a}}|\mathbf{a}^\diamond|\mathbf{c}_s)$ for some $s$, or $\mathbf{p}_i$ is a knot appearing in $\mathfrak{r}_u(\mathbf{k},\mathbf{a})$ and $c_i\in R$ is its coefficient.
In the same sprit, we will rewrite (\ref{eq:substitution-2}) as $\varepsilon_d(\mathbf{k},\mathbf{a})=\sum_j\mathbf{q}_jd_j$.
Always remember that $\mathbf{p}_i,\mathbf{q}_j\prec\mathbf{k}$ for all $i,j$.
\end{rmk}

\begin{figure}[H]
  \centering
  \includegraphics[width=12.4cm]{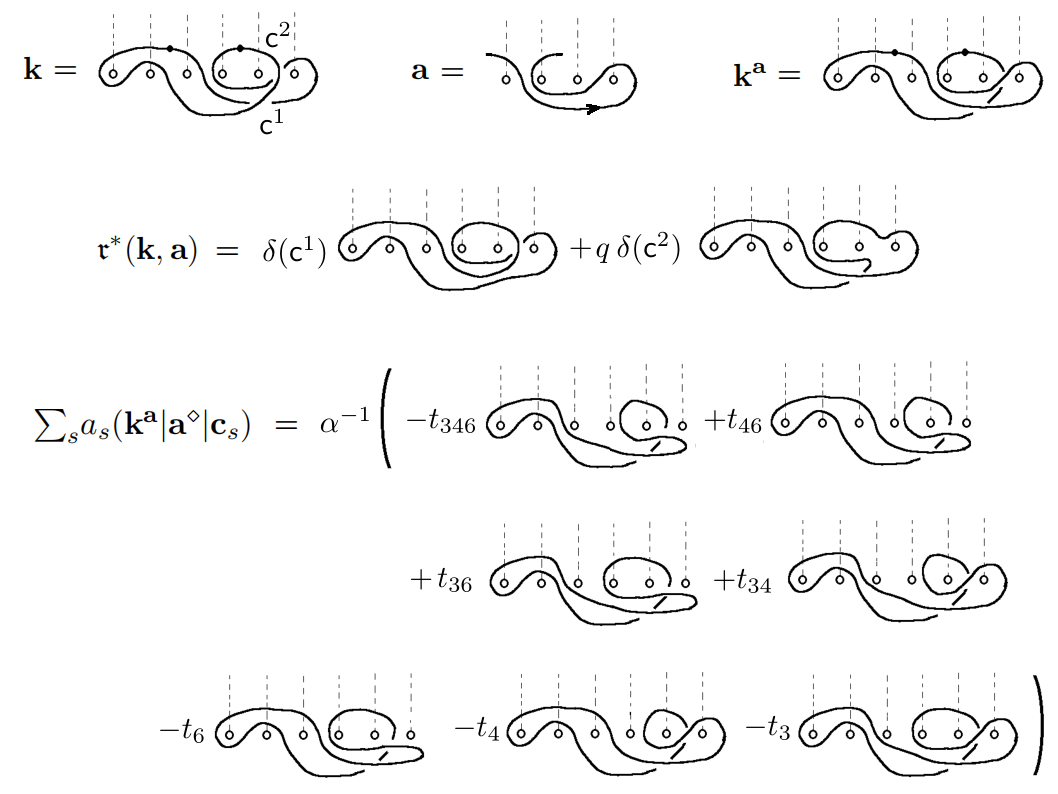}
  \caption{Note that $\epsilon(\mathsf{c}^1)=1$ and $\epsilon(\mathsf{c}^2)=-1$. The third to fifth rows present ${\sum}_sa_s(\mathbf{k}^{\mathbf{a}}|\mathbf{a}^\diamond|\mathbf{c}_s)$, which is obtained by replacing
  $\mathbf{a}^\diamond$ by ${\rm ch}_u(\mathbf{a}^\diamond)$.}\label{fig:example-0}
\end{figure}

\begin{exmp}\label{exmp:subsititute}
\rm A degree $7$ knot $\mathbf{k}$ is given in the upper-left corner of Figure \ref{fig:example-0}. Let $\mathbf{a}$ denote the degree $3$ arc bounded by the two bullets. An orientation of $\mathbf{a}$ is chosen, so that elements of ${\rm Cr}(\langle\mathbf{k}|\mathbf{a}\rangle,\mathbf{a})$ can be enumerated as $\mathsf{c}^1,\mathsf{c}^2$.
In the case $\alpha^{-1}\in R$, applying the formula for ${\rm ch}_u(\mathbf{a}^\diamond)$ given in Figure \ref{fig:curve}, with a representative $\mathbf{c}_s$ chosen for each $[\mathbf{c}_s]$, we obtain ${\sum}_sa_s(\mathbf{k}^{\mathbf{a}}|\mathbf{a}^\diamond|\mathbf{c}_s)$.
\end{exmp}

\begin{lem}
Suppose $\mathbf{k}$ is a knot.
Then $\varepsilon_u(\mathbf{k},\mathbf{a})\equiv\varepsilon_d(\mathbf{k},\mathbf{b})\pmod{\mathcal{L}_{\mathbf{k}}}$
for any $\mathbf{a},\mathbf{b}\in{\rm Ar}^\ast(\mathbf{k})$.
\end{lem}

\begin{proof}
The assertion holds automatically when $|\mathbf{k}|\le 2h$.

We suppose $|\mathbf{k}|>2h$. Abbreviate $g\equiv g'\pmod{\mathcal{L}_{\mathbf{k}}}$ to $g\equiv g'$.

{\bf Step 1.} Suppose $\mathbf{a},\mathbf{b}\in{\rm Ar}^\ast(\mathbf{k})$ with $\mathbf{a}\cap\mathbf{b}=\emptyset$.

\begin{figure}[h]
  \centering
  \includegraphics[width=11cm]{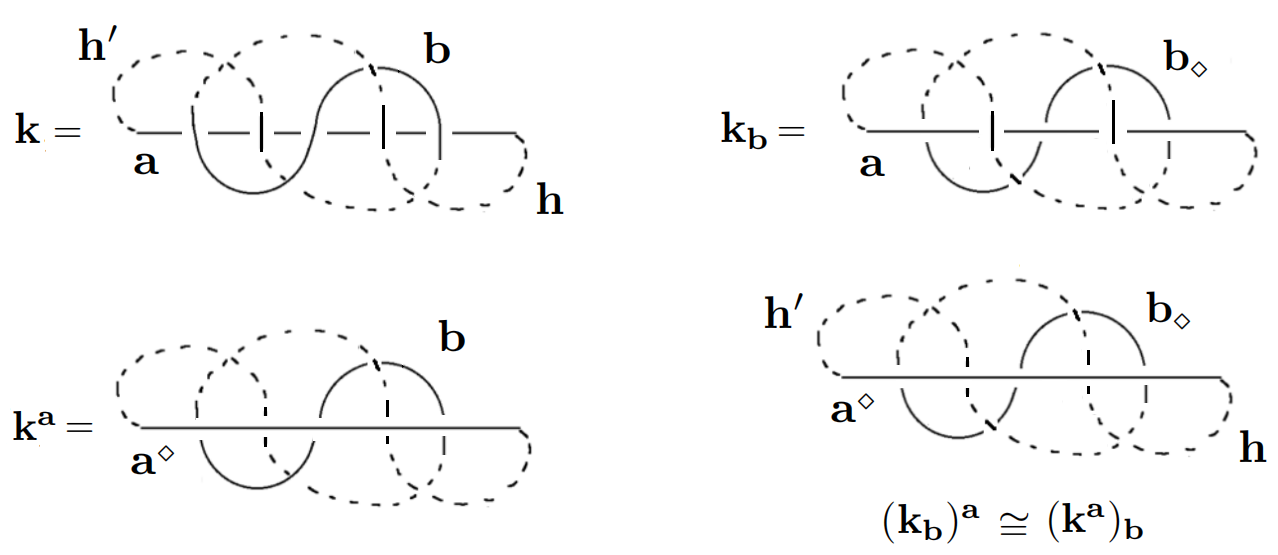}\\
  \caption{Shown at the upper-left corner is a knot $\mathbf{k}$, in which the horizontal line presents $\mathbf{a}$, the solid arc presents $\mathbf{b}$, and the two dotted arcs stand for $\mathbf{h},\mathbf{h}'$. Shown at the lower-right corner is
  $(\mathbf{k}_{\mathbf{b}})^{\mathbf{a}}\cong(\mathbf{k}^{\mathbf{a}})_{\mathbf{b}}\cong{\rm tr}(\mathbf{a}^\diamond \mathbf{h}\mathbf{b}_\diamond \mathbf{h}')$.}\label{fig:explain}
\end{figure}

Write $\mathbf{k}={\rm tr}(\mathbf{a}\mathbf{h}\mathbf{b}\mathbf{h}')$, for some arcs $\mathbf{h},\mathbf{h}'$.
Denote the arc of $\mathbf{k}_{\mathbf{b}}$ inherited from $\mathbf{a}$ also by $\mathbf{a}$, and denote the arc of $\mathbf{k}^{\mathbf{a}}$ inherited from $\mathbf{b}$ also by $\mathbf{b}$. We will adopt such conventions in similar situations.
Clearly, $(\mathbf{k}_{\mathbf{b}})^{\mathbf{a}}\cong(\mathbf{k}^{\mathbf{a}})_{\mathbf{b}}
\cong{\rm tr}(\mathbf{a}^\diamond \mathbf{h}\mathbf{b}_\diamond \mathbf{h}')$.
See Figure \ref{fig:explain} for illustration.

Suppose ${\rm ch}_u(\mathbf{a}^\diamond)=\sum_sa_s[\mathbf{c}_s]$ and ${\rm ch}_d(\mathbf{b}_\diamond)=\sum_t[\mathbf{d}_t]b_t$.
Let
$$\mathbf{j}_{st}={\rm tr}(\mathbf{c}_s\mathbf{h}\mathbf{d}_t\mathbf{h}')={\rm tr}(\mathbf{d}_t\mathbf{h}'\mathbf{c}_s\mathbf{h}).$$
Let $\hat{\epsilon}=\hat{\epsilon}(\mathbf{k},\mathbf{a})$, $\hat{\epsilon}'=\hat{\epsilon}(\mathbf{k}_{\mathbf{b}},\mathbf{a})$, $\check{\epsilon}=\check{\epsilon}(\mathbf{k},\mathbf{b})$, $\check{\epsilon}'=\check{\epsilon}(\mathbf{k}^{\mathbf{a}},\mathbf{b})$.
Then
\begin{align*}
\hat{\epsilon}+\check{\epsilon}'=\hat{\epsilon}+\check{\epsilon}
-{\sum}_{\mathsf{c}\in{\rm Cr}(\mathbf{b},\mathbf{a})}\epsilon(\mathsf{c})=\check{\epsilon}+\hat{\epsilon}'.
\end{align*}

For each $s$, let $\mathbf{r}_s=(\mathbf{k}^{\mathbf{a}}|\mathbf{a}^\diamond|\mathbf{c}_s)$;
since $(\mathbf{r}_s)_{\mathbf{b}}={\rm tr}(\mathbf{c}_s\mathbf{h}\mathbf{b}_\diamond\mathbf{h}')$,
we have
$((\mathbf{r}_s)_{\mathbf{b}}|\mathbf{b}_\diamond|\mathbf{d}_t)={\rm tr}(\mathbf{c}_s\mathbf{h}\mathbf{d}_t\mathbf{h}')=\mathbf{j}_{st}$,
hence
$$\mathbf{r}_s=\mathfrak{r}_d(\mathbf{r}_s,\mathbf{b})+q^{\check{\epsilon}'}{\sum}_t\mathbf{j}_{st}b_t \quad \text{in\ \ }\mathcal{S}_n.$$
The knots appearing in $\mathfrak{r}_d(\mathbf{k}^{\mathbf{a}},\mathbf{b})$ bijectively correspond to those in $\mathfrak{r}_d(\mathbf{r}_s,\mathbf{b})$ via $\mathbf{k}'\mapsto(\mathbf{k}'|\mathbf{a}^\diamond|\mathbf{c}_s)$, with the same coefficients.
Each such knot $\mathbf{k}'$ satisfies $\mathbf{k}'\prec\mathbf{k}$ and $\mathbf{k}'=\sum_sa_s(\mathbf{k}'|\mathbf{a}^\diamond|\mathbf{c}_s)$ in $\mathcal{S}_n$, so $\mathbf{k}'-\sum_sa_s(\mathbf{k}'|\mathbf{a}^\diamond|\mathbf{c}_s)$ is a type (ii) element of $\mathcal{L}_{\mathbf{k}}$, implying $\mathbf{k}'\equiv\sum_sa_s(\mathbf{k}'|\mathbf{a}^\diamond|\mathbf{c}_s)$.
Hence
$$\mathfrak{r}_d(\mathbf{k}^{\mathbf{a}},\mathbf{b})\equiv{\sum}_sa_s\mathfrak{r}_d(\mathbf{r}_s,\mathbf{b}).$$
Consequently,
\begin{align*}
{\sum}_sa_s\mathbf{r}_s\equiv {\sum}_sa_s\left(\mathfrak{r}_d(\mathbf{r}_s,\mathbf{b})+q^{\check{\epsilon}'}{\sum}_t\mathbf{j}_{st}b_t\right)
\equiv\mathfrak{r}_d(\mathbf{k}^{\mathbf{a}},\mathbf{b})+q^{\check{\epsilon}'}{\sum}_{s,t}a_s\mathbf{j}_{st}b_t.
\end{align*}
Recalling $\varepsilon_u(\mathbf{k},\mathbf{a})=\mathfrak{r}_u(\mathbf{k},\mathbf{a})+q^{\hat{\epsilon}}{\sum}_sa_s\mathbf{r}_s$, we obtain
\begin{align*}
\varepsilon_u(\mathbf{k},\mathbf{a})
\equiv\mathfrak{r}_u(\mathbf{k},\mathbf{a})+q^{\hat{\epsilon}}\mathfrak{r}_d(\mathbf{k}^{\mathbf{a}},\mathbf{b})
+q^{\hat{\epsilon}+\check{\epsilon}'}{\sum}_{s,t}a_s\mathbf{j}_{st}b_t.
\end{align*}

In a parallel way, we can deduce
$$\varepsilon_d(\mathbf{k},\mathbf{b})
\equiv\mathfrak{r}_d(\mathbf{k},\mathbf{b})+q^{\check{\epsilon}}\mathfrak{r}_u(\mathbf{k}_{\mathbf{b}},\mathbf{a})
+q^{\check{\epsilon}+\hat{\epsilon}'}{\sum}_{s,t}a_s\mathbf{j}_{st}b_t.$$

Thus,
$\varepsilon_u(\mathbf{k},\mathbf{a})-\varepsilon_d(\mathbf{k},\mathbf{b})\equiv g$, with
$$g=\mathfrak{r}_u(\mathbf{k},\mathbf{a})+q^{\hat{\epsilon}}\mathfrak{r}_d(\mathbf{k}^{\mathbf{a}},\mathbf{b})
-\mathfrak{r}_d(\mathbf{k},\mathbf{b})-q^{\check{\epsilon}}\mathfrak{r}_u(\mathbf{k}_{\mathbf{b}},\mathbf{a}).$$
Since $g$ is a $R$-linear combination of knots simpler than $\mathbf{k}$ and
$$\tilde{\theta}(g)=\tilde{\theta}(\varepsilon_u(\mathbf{k},\mathbf{b}))-\tilde{\theta}(\varepsilon_d(\mathbf{k},\mathbf{b}))
=[\mathbf{k}]-[\mathbf{k}]=0,$$
we have $g\in\mathcal{L}_{\mathbf{k}}$, so that $\varepsilon_u(\mathbf{k},\mathbf{a})\equiv\varepsilon_d(\mathbf{k},\mathbf{b})$.

\medskip

{\bf Step 2.} Suppose $\mathbf{a}_1,\mathbf{a}_2\in{\rm Ar}^\ast(\mathbf{k})$ are {\it successive}, by which we mean
$$\max\{|\mathbf{a}_1|,|\mathbf{a}_2|\}<|\mathbf{a}_1\cup\mathbf{a}_2|\le h+1.$$
Due to $|\mathbf{k}|>2h$, there exists
$\mathbf{b}\in{\rm Ar}^\ast(\mathbf{k})$ with $\mathbf{b}\cap\mathbf{a}_1=\mathbf{b}\cap\mathbf{a}_2=\emptyset$.
Hence
\begin{align*}
\varepsilon_u(\mathbf{k},\mathbf{a}_1)&\equiv\varepsilon_d(\mathbf{k},\mathbf{b})\equiv\varepsilon_u(\mathbf{k},\mathbf{a}_2),   \\
\varepsilon_d(\mathbf{k},\mathbf{a}_1)&\equiv\varepsilon_u(\mathbf{k},\mathbf{b})\equiv\varepsilon_d(\mathbf{k},\mathbf{a}_2).
\end{align*}

On the other hand, for each $\mathbf{a}\in{\rm Ar}^\ast(\mathbf{k})$, take $\mathbf{c}\in{\rm Ar}(\mathbf{k})$ such that $\mathbf{c}\cap\mathbf{a}$ consists of a single point and $|\mathbf{c}\cup\mathbf{a}|=h+1$, then there exists $\mathbf{a}'\in{\rm Ar}^\ast(\mathbf{k})$ such that $\mathbf{a}'\subseteq\mathbf{c}\cup\mathbf{a}$ and $\mathbf{a},\mathbf{a}'$ are successive.
Thus, any $\mathbf{a}_1,\mathbf{a}_2\in{\rm Ar}^\ast(\mathbf{k})$ can be related by a string of successive pairs,
so that $\varepsilon_u(\mathbf{k},\mathbf{a}_1)\equiv\varepsilon_u(\mathbf{k},\mathbf{a}_2)$, and
$\varepsilon_d(\mathbf{k},\mathbf{a}_1)\equiv\varepsilon_d(\mathbf{k},\mathbf{a}_2)$.

For any $\mathbf{a}\in{\rm Ar}^\ast(\mathbf{k})$, take $\mathbf{b}\in{\rm Ar}^\ast(\mathbf{k})$ with $\mathbf{a}\cap \mathbf{b}=\emptyset$. Then by the above results, $\varepsilon_u(\mathbf{k},\mathbf{a})\equiv\varepsilon_d(\mathbf{k},\mathbf{b})\equiv\varepsilon_d(\mathbf{k},\mathbf{a})$.

The proof is completed.
\end{proof}

\begin{defn}
\rm For a knot $\mathbf{k}$, define $\varepsilon(\mathbf{k})$ to be the common value of $\varepsilon_u(\mathbf{k},\mathbf{a})\pmod{\mathcal{L}_{\mathbf{k}}}$ and $\varepsilon_d(\mathbf{k},\mathbf{a})\pmod{\mathcal{L}_{\mathbf{k}}}$ for an arbitrary $\mathbf{a}\in{\rm Ar}^\ast(\mathbf{k})$.

For a stacked link $\mathbf{l}=\mathbf{k}_1\cdots\mathbf{k}_r$ with $r\ge 2$, define $\varepsilon(\mathbf{l})$ as $\varepsilon(\mathbf{k}_1)\cdots\varepsilon(\mathbf{k}_r)\pmod{\mathcal{L}_{|\mathbf{l}|,0}}$.
\end{defn}

\begin{rmk}\label{rmk:chop-up}
\rm (i) We call $\varepsilon$ the {\it chopping-up} map.
Intuitively, $\varepsilon(\mathbf{l})$ is the result of chopping up an arbitrary arc in ${\rm Ar}^\ast(\mathbf{l})$ for each component of $\mathbf{l}$.

It is clear that $\varepsilon(\mathbf{k})\equiv\varepsilon(\mathbf{k'})\pmod{\mathcal{L}_{\mathbf{k}}}$ if the knots $\mathbf{k}\cong\mathbf{k}'$.
Also clear is that $\varepsilon(\mathbf{l})=\mathbf{l}$ if $\mathbf{l}\in\mathcal{T}$.

(ii) From the definitions (in particular, the type (ii) and (iii) elements in $\mathcal{L}_{m,c}$) we see that $\mathbf{k}\equiv\varepsilon(\mathbf{k})\pmod{\mathcal{L}_{m,c}}$ for any knot $\mathbf{k}$ with $\|\mathbf{k}\|\prec(m,c)$.
\end{rmk}

\section{Properties of the chopping-up map}

\begin{lem}\label{lem:crossing-change}
Suppose $\mathbf{k}$ is a knot.
Then for each $\mathsf{c}\in{\rm Cr}(\mathbf{k})$,
$$\varepsilon(\mathbf{k})\equiv q^{\epsilon(\mathsf{c})}\varepsilon(\mathbf{k}^\star(\mathsf{c}))
+\delta(\mathsf{c})\varepsilon(\mathbf{k}^\flat(\mathsf{c}))\pmod{\mathcal{L}_{\mathbf{k}}}.$$
Consequently, for each $\mathbf{a}\in{\rm Ar}(\mathbf{k})$,
$$\varepsilon(\mathbf{k})\equiv
q^{\hat{\epsilon}(\mathbf{k},\mathbf{a})}\varepsilon(\mathbf{k}^{\mathbf{a}})+\varepsilon(\mathfrak{r}_u(\mathbf{k},\mathbf{a}))
\equiv q^{\check{\epsilon}(\mathbf{k},\mathbf{a})}\varepsilon(\mathbf{k}_{\mathbf{a}})+\varepsilon(\mathfrak{r}_d(\mathbf{k},\mathbf{a}))
\pmod{\mathcal{L}_{\mathbf{k}}},$$
where $\varepsilon(\mathfrak{r}_u(\mathbf{k},\mathbf{a}))$, $\varepsilon(\mathfrak{r}_d(\mathbf{k},\mathbf{a}))$ are defined by $R$-linear extensions.
\end{lem}

\begin{proof}
The assertion holds automatically if $|\mathbf{k}|\le 2h$, due to the definition of $\mathcal{L}_{m,c}$ for $m\le 2h$.
Indeed, $\varepsilon(\mathbf{k})-q^{\epsilon(\mathsf{c})}\varepsilon(\mathbf{k}^\star(\mathsf{c}))
-\delta(\mathsf{c})\varepsilon(\mathbf{k}^\flat(\mathsf{c}))=0$ in $\mathcal{S}_n$, and each stacked link appearing in
$\varepsilon(\mathbf{k})$, $\varepsilon(\mathbf{k}^\star(\mathsf{c}))$, $\varepsilon(\mathbf{k}^\flat(\mathsf{c}))$ has degree $\le 2h$.

Suppose $|\mathbf{k}|>2h$.
Take $\mathbf{a}\in{\rm Ar}^\ast(\mathbf{k})$ away from $\mathsf{c}$, by which we mean
${\rm over}(\mathsf{c}),{\rm under}(\mathsf{c})\notin\mathbf{a}$.
Write $\varepsilon_u(\mathbf{k},\mathbf{a})\equiv\sum_ic_i\mathbf{p}_i\pmod{\mathcal{L}_{\mathbf{k}}}$ as in Remark \ref{rmk:convention} (iii), then
$\mathbf{p}_i,\mathbf{p}_i^\star(\mathsf{c}),\mathbf{p}_i^\flat(\mathsf{c})\prec\mathbf{k}$ for each $i$, and
\begin{align*}
\varepsilon_u(\mathbf{k}^\star(\mathsf{c}),\mathbf{a})&\equiv{\sum}_ic_i\mathbf{p}_i^\star(\mathsf{c})\pmod{\mathcal{L}_{\mathbf{k}}},   \\
\varepsilon_u(\mathbf{k}^\flat(\mathsf{c}),\mathbf{a})&\equiv{\sum}_ic_i\mathbf{p}_i^\flat(\mathsf{c})\pmod{\mathcal{L}_{\mathbf{k}}}.
\end{align*}
Since $\mathbf{p}_i=q^{\epsilon(\mathsf{c})}\mathbf{p}_i^\star(\mathsf{c})+\delta(\mathsf{c})\mathbf{p}_i^\flat(\mathsf{c})$ in $\mathcal{S}_n$ for each $i$, we have
$$\varepsilon_u(\mathbf{k},\mathbf{a})\equiv q^{\epsilon(\mathsf{c})}\varepsilon_u(\mathbf{k}^\star(\mathsf{c}),\mathbf{a})
+\delta(\mathsf{c})\cdot\varepsilon_u(\mathbf{k}^\flat(\mathsf{c}),\mathbf{a})\pmod{\mathcal{L}_{\mathbf{k}}}.$$
Thus, $\varepsilon(\mathbf{k})\equiv q^{\epsilon(\mathsf{c})}\varepsilon(\mathbf{k}^\star(\mathsf{c}))
+\delta(\mathsf{c})\varepsilon(\mathbf{k}^\flat(\mathsf{c}))\pmod{\mathcal{L}_{\mathbf{k}}}.$
\end{proof}

\begin{defn}\label{defn:EST}
\rm Call a skein triple $(\mathbf{l}_\times,\mathbf{l}_\infty,\mathbf{l}_0)$ an {\it EST}, if one of $\mathbf{l}_\times,\mathbf{l}_\infty,\mathbf{l}_0$ is a two-component stacked link and the other two are knots.
\end{defn}


\begin{lem}\label{lem:EST}
For each EST $(\mathbf{l}_\times,\mathbf{l}_\infty,\mathbf{l}_0)$,
$$\varepsilon(\mathbf{l}_\times)\equiv q^{\frac{1}{2}}\varepsilon(\mathbf{l}_\infty)+\overline{q}^{\frac{1}{2}}\varepsilon(\mathbf{l}_0)
\pmod{\mathcal{L}_{\mathbf{l}_\times}}.$$
\end{lem}

\begin{proof}
The assertion holds automatically if $|\mathbf{l}_\times|\le 2h$. Suppose $|\mathbf{l}_\times|>2h$.

Suppose $\mathbf{l}_\times=\mathbf{k}_1\mathbf{k}_2$, $\mathbf{l}_\infty=(\mathbf{k}_1\mathbf{k}_2)^{\mathsf{c}}_\infty$ and $\mathbf{l}_0=(\mathbf{k}_1\mathbf{k}_2)^{\mathsf{c}}_0$ for some $\mathsf{c}\in{\rm Cr}(\mathbf{k}_1,\mathbf{k}_2)$;
the other two cases are similar.

Since $|\mathbf{l}_\times|>2h$, we have $|\mathbf{k}_1|>h$ or $|\mathbf{k}_2|>h$. Assume $|\mathbf{k}_1|>h$; the other case is similar.
Use induction on ${\rm cn}(\mathbf{k}_1)$ to prove the assertion.

For $g,g'\in\mathcal{F}$, we abbreviate $g\equiv g'\pmod{\mathcal{L}_{\mathbf{l}_\times}}$ to $g\equiv g'$.

Take $\mathbf{a}\in{\rm Ar}^\ast(\mathbf{k}_1)$ away from $\mathsf{c}$. Write $\varepsilon_u(\mathbf{k}_1,\mathbf{a})\equiv\sum_ic_i\mathbf{p}_i\pmod{\mathcal{L}_{\mathbf{k}_1}}$ as in Remark \ref{rmk:convention} (iii). Let $\mathbf{u}_i=\mathbf{p}_i\mathbf{k}_2$. Then
$$\varepsilon(\mathbf{l}_\times)\equiv{\sum}_ic_i\mathbf{u}_i, \qquad \varepsilon(\mathbf{l}_\infty)\equiv{\sum}_ic_i\cdot(\mathbf{u}_i)^\mathsf{c}_\infty, \qquad \varepsilon(\mathbf{l}_0)\equiv{\sum}_ic_i\cdot(\mathbf{u}_i)^\mathsf{c}_0;$$
the second and third equations are due to that the operations applied to $\mathbf{a}$ do not interfere resolving $\mathsf{c}$.

If ${\rm cn}(\mathbf{k}_1)=0$, then for each $i$, we have $c_i\in\mathcal{T}$ and $|\mathbf{p}_i|<|\mathbf{k}_1|$, so $$c_i\big(\mathbf{u}_i-q^{\frac{1}{2}}(\mathbf{u}_i)^\mathsf{c}_\infty-\overline{q}^{\frac{1}{2}}(\mathbf{u}_i)^\mathsf{c}_0\big)
\in\mathcal{L}_{\mathbf{l}_\times},$$
implying
$$\varepsilon(\mathbf{l}_\times)\equiv{\sum}_ic_i\left(q^{\frac{1}{2}}(\mathbf{u}_i)^\mathsf{c}_\infty
+\overline{q}^{\frac{1}{2}}(\mathbf{u}_i)^\mathsf{c}_0\right)\equiv q^{\frac{1}{2}}\varepsilon(\mathbf{l}_\infty)+\overline{q}^{\frac{1}{2}}\varepsilon(\mathbf{l}_0).$$
Suppose ${\rm cn}(\mathbf{k}_1)>0$ and the assertion is true when ${\rm cn}(\mathbf{k}_1)$ is smaller. There probably exists $i$ such that $c_i\in R$, $|\mathbf{p}_i|=|\mathbf{k}_1|$ and ${\rm cn}(\mathbf{p}_i)<{\rm cn}(\mathbf{k}_1)$.
Then for the EST
$(\mathbf{u}_i,(\mathbf{u}_i)^\mathsf{c}_\infty,(\mathbf{u}_i)^\mathsf{c}_0)$, by the inductive hypothesis,
$$\varepsilon(\mathbf{u}_i)\equiv q^{\frac{1}{2}}\varepsilon((\mathbf{u}_i)^\mathsf{c}_\infty)
+\overline{q}^{\frac{1}{2}}\varepsilon((\mathbf{u}_i)^\mathsf{c}_0).$$
On the other hand,
$$\mathbf{u}_i=\mathbf{p}_i\mathbf{k}_2\equiv\varepsilon(\mathbf{p}_i)\mathbf{k}_2\equiv\varepsilon(\mathbf{u}_i),$$
and by Remark \ref{rmk:chop-up} (ii), $\mathbf{j}\equiv\varepsilon(\mathbf{j})$ for $\mathbf{j}\in\{(\mathbf{u}_i)^\mathsf{c}_\infty,(\mathbf{u}_i)^\mathsf{c}_0\}$. Hence
$$\mathbf{u}_i\equiv q^{\frac{1}{2}}(\mathbf{u}_i)^\mathsf{c}_\infty+\overline{q}^{\frac{1}{2}}(\mathbf{u}_i)^\mathsf{c}_0.$$
Thus, still $\varepsilon(\mathbf{l}_\times)\equiv
q^{\frac{1}{2}}\varepsilon(\mathbf{l}_\infty)+\overline{q}^{\frac{1}{2}}\varepsilon(\mathbf{l}_0)$.
\end{proof}

\begin{lem}\label{lem:disjoint}
$\varepsilon(\mathbf{k}_1)\varepsilon(\mathbf{k}_2)\equiv\varepsilon(\mathbf{k}_2)\varepsilon(\mathbf{k}_1)\pmod{\mathcal{L}_{m,c}}$ for any knots $\mathbf{k}_1,\mathbf{k}_2$ with $\pi(\mathbf{k}_1)\cap\pi(\mathbf{k}_2)=\emptyset$ and $\|\mathbf{k}_1\mathbf{k}_2\|\prec(m,c)$.
\end{lem}

\begin{proof}
The projection $\pi(\mathbf{k}_1\cup \mathbf{k}_2)$ cuts $\Sigma\setminus\mathbf{z}$ into disks, among which there must be at least one, denoted by $\Upsilon$, satisfying $\partial \Upsilon\cap\pi(\mathbf{k}_i)\ne\emptyset$ for $i=1,2$.
Take a short $\mathbf{c}_i\in{\rm Ar}(\mathbf{k}_i)$ with $\pi(\mathbf{c}_i)\subset\partial \Upsilon$, and let $\mathbf{k}_i^\circ=\mathbf{k}_i\setminus\mathring{\mathbf{c}_i}$.

\begin{figure}[h]
  \centering
  \includegraphics[width=13cm]{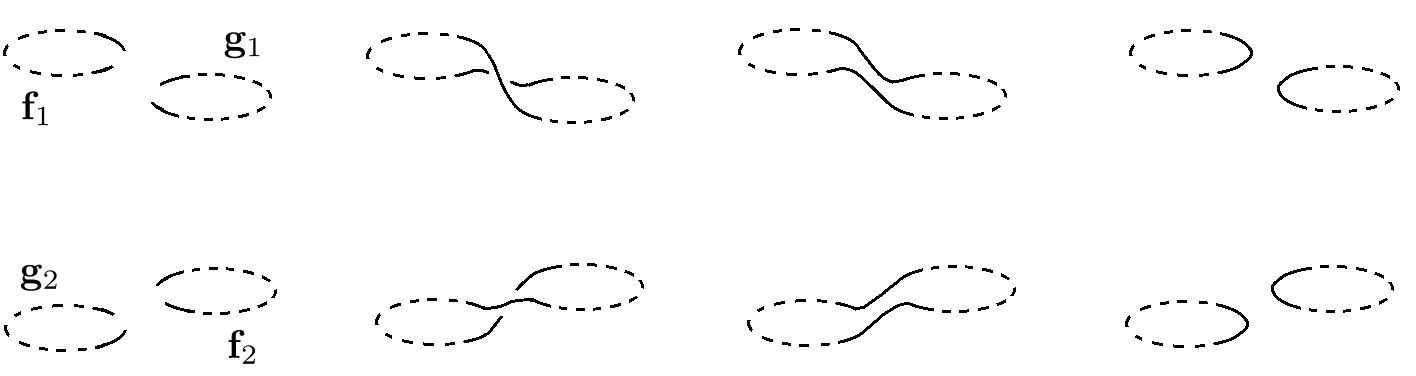}\\
  \caption{First row (from left to right): $\mathbf{f}_1,\mathbf{g}_1$; $\mathbf{t}_1$;
  $\mathbf{t}'_1$; $\mathbf{k}_1\mathbf{k}_2$.
  Second row (from left to right): $\mathbf{g}_2,\mathbf{f}_2$; $\mathbf{t}_2$; $\mathbf{t}'_2$; $\mathbf{k}_2\mathbf{k}_1$. The dotted arcs stand for the remaining parts which are irrelevant. We arrange that $\mathbf{f}_1,\mathbf{f}_2\subset\Sigma\times(\frac{1}{2},1)$, and $\mathbf{g}_1,\mathbf{g}_2\subset\Sigma\times(0,\frac{1}{2})$.}
  \label{fig:A-B}
\end{figure}

Let $\mathbf{f}_1\subset\Sigma\times(\frac{1}{2},1)$ be a copy of $\mathbf{k}_1^\circ$, by which we mean the image of $\mathbf{k}_1^\circ$ under the evident map $\Sigma\times(0,1)\to\Sigma\times(\frac{1}{2},1)$, and let $\mathbf{g}_1\subset\Sigma\times(0,\frac{1}{2})$ be a copy of $\mathbf{k}_2^\circ$. Construct $\mathbf{t}_1$, $\mathbf{t}'_1$ using simple arcs in $\Upsilon\times(0,1)$ as illustrated in the first row of Figure \ref{fig:A-B}. Then
$(\mathbf{t}_1,\mathbf{t}'_1,\mathbf{k}_1\mathbf{k}_2)$
is an EST. By Lemma \ref{lem:EST},
\begin{align*}
\varepsilon(\mathbf{t}_1)\equiv q^{\frac{1}{2}}\varepsilon(\mathbf{t}'_1)+\overline{q}^{\frac{1}{2}}\varepsilon(\mathbf{k}_1\mathbf{k}_2)
\pmod{\mathcal{L}_{m,c}}.
\end{align*}

Let $\mathbf{g}_2\subset\Sigma\times(0,\frac{1}{2})$ be a copy of $\mathbf{k}_1^\circ$, and $\mathbf{f}_2\subset\Sigma\times(\frac{1}{2},1)$ a copy of $\mathbf{k}_2^\circ$. Construct $\mathbf{t}_2$, $\mathbf{t}'_2$ using simple arcs in $\Upsilon\times(0,1)$ as illustrated in the second row of Figure \ref{fig:A-B}.
Then $(\mathbf{t}_2,\mathbf{t}'_2,\mathbf{k}_2\mathbf{k}_1)$ is an EST.
By Lemma \ref{lem:EST},
\begin{align*}
\varepsilon(\mathbf{t}_2)\equiv q^{\frac{1}{2}}\varepsilon(\mathbf{t}'_2)+\overline{q}^{\frac{1}{2}}\varepsilon(\mathbf{k}_2\mathbf{k}_1)
\pmod{\mathcal{L}_{m,c}}.
\end{align*}

Note that $\mathbf{t}_1\cong \mathbf{t}_2$, $\mathbf{t}'_1\cong \mathbf{t}'_2$,
so $\varepsilon(\mathbf{t}_1)\equiv\varepsilon(\mathbf{t}_2)\pmod{\mathcal{L}_{m,c}}$, and
$\varepsilon(\mathbf{t}'_1)\equiv\varepsilon(\mathbf{t}'_2)\pmod{\mathcal{L}_{m,c}}$.
Therefore, $\varepsilon(\mathbf{k}_1\mathbf{k}_2)\equiv\varepsilon(\mathbf{k}_2\mathbf{k}_1)\pmod{\mathcal{L}_{m,c}}$.
\end{proof}


\begin{cor}\label{cor:multicurve}
If $\mathbf{m}=\mathbf{s}_1\sqcup\cdots\sqcup \mathbf{s}_r$ is a multicurve with $|\mathbf{m}|\le m$, then
$\varepsilon(\mathbf{s}_{\sigma(1)}\cdots \mathbf{s}_{\sigma(r)})\equiv\varepsilon(\mathbf{s}_1\cdots\mathbf{s}_r)\pmod{\mathcal{L}_{m,1}}$ for each $\sigma\in{\rm Sym}(r)$.
\end{cor}

\begin{proof}
For each $k<r$, by Lemma \ref{lem:disjoint},
$\varepsilon(\mathbf{s}_{k}\mathbf{s}_{k+1})\equiv\varepsilon(\mathbf{s}_{k+1}\mathbf{s}_{k})\pmod{\mathcal{L}_{m,1}}$,
hence
$$\varepsilon(\mathbf{s}_{1})\cdots\varepsilon(\mathbf{s}_{k})\varepsilon(\mathbf{s}_{k+1})\cdots\varepsilon(\mathbf{s}_{r})\equiv
\varepsilon(\mathbf{s}_{1})\cdots\varepsilon(\mathbf{s}_{k+1})\varepsilon(\mathbf{s}_{k})\cdots\varepsilon(\mathbf{s}_{r})
\pmod{\mathcal{L}_{m,1}}.$$
The assertion follows from the basic fact that each $\sigma\in{\rm Sym}(r)$ is the composite of transpositions of the form $(k,k+1)$.
\end{proof}

\begin{lem}\label{lem:simple}
If $\mathbf{s},\mathbf{s}'\subset\Sigma$ are simple curves such that $|\mathbf{s}|\le m$ and $\mathbf{s}'$ results from reducing a shrinkable arc of $\mathbf{s}$, then $\varepsilon(\mathbf{s})\equiv\varepsilon(\mathbf{s}')\pmod{\mathcal{L}_{m,0}}$.
\end{lem}

\begin{proof}
The assertion automatically holds when $|\mathbf{s}|\le 2h$.

Suppose $|\mathbf{s}|>2h$. Let $\mathbf{b}\in{\rm Ar}_2(\mathbf{s})$ be shrinkable.
Take $\mathbf{a}\in{\rm Ar}^\ast(\langle\mathbf{s}|\mathbf{b}\rangle)$, and write ${\rm ch}_u(\mathbf{a})=\sum_sa_s[\mathbf{c}_s]$. Then
$$\varepsilon_u(\mathbf{s},\mathbf{a})={\sum}_sa_s\cdot(\mathbf{s}|\mathbf{a}|\mathbf{c}_s)
\equiv{\sum}_sa_s\cdot(\mathbf{s}'|\mathbf{a}|\mathbf{c}_s)=\varepsilon_u(\mathbf{s}',\mathbf{a})\pmod{\mathcal{L}_{m,0}},$$
as $(\mathbf{s}|\mathbf{a}|\mathbf{c}_s)=(\mathbf{s}'|\mathbf{a}|\mathbf{c}_s)$ in $\mathcal{S}_n$ for each $s$.
Hence $\varepsilon(\mathbf{s})\equiv\varepsilon(\mathbf{s}')\pmod{\mathcal{L}_{m,0}}$.
\end{proof}

\begin{rmk}\label{rmk:multicurve}
\rm The upshot is: for a multicurve $\mathbf{m}=\mathbf{s}_1\sqcup\cdots\sqcup\mathbf{s}_r$ with $|\mathbf{m}|\le m$, we can define $\varepsilon(\mathbf{m})$ to be $\varepsilon(\mathbf{s}_1\cdots\mathbf{s}_r)\pmod{\mathcal{L}_{m,1}}$.

Modulo $\mathcal{L}_{m,1}$, the notion $\varepsilon(g)$ can be extended (by $R$-linearity), to allow $g$ to be a $R$-linear combination of multicurves of degree
$\le m$.

Recall that $\Theta(\mathbf{m})$ is the irreducible multicurve obtained by reducing shrinkable arcs and removing nullhomotopic circles in $\mathbf{m}$. By Lemma \ref{lem:simple}, $\varepsilon(\mathbf{m})\equiv\varepsilon(\Theta(\mathbf{m}))\pmod{\mathcal{L}_{m,1}}$.

If $\mathbf{m}_1,\mathbf{m}_2$ are isotopic multicurves of degree $\le m$, then by Lemma \ref{lem:reduce-multi-curve},
$\Theta(\mathbf{m}_1)\cong\Theta(\mathbf{m}_2)$, hence $\varepsilon(\mathbf{m}_1)\equiv\varepsilon(\mathbf{m}_2)\pmod{\mathcal{L}_{m,1}}$.

Consequently, if $g_1,g_2$ are $R$-linear combinations of (possibly reducible) multicurves with $g_1=g_2$ in $\mathcal{S}_n$, then $\varepsilon(g_1)\equiv\varepsilon(g_2)\pmod{\mathcal{L}_{m,1}}$, where $m$ is the maximum of the degrees of the multicurves appearing in $g_1,g_2$.
\end{rmk}

\begin{lem}\label{lem:resolve}
$\varepsilon(\mathbf{k})\equiv\varepsilon(\Theta(\mathbf{k}))\pmod{\mathcal{L}_{\mathbf{k}}}$ for each knot $\mathbf{k}$.
\end{lem}

\begin{proof}
When $\mathbf{k}$ is simple, trivially $\varepsilon(\mathbf{k})\equiv\varepsilon(\Theta(\mathbf{k}))\pmod{\mathcal{L}_{\mathbf{k}}}$.

Suppose ${\rm cn}(\mathbf{k})>0$. We first show that $$\varepsilon(\mathbf{k})\equiv{\sum}_jb_{j}\varepsilon(\mathbf{g}_{j})\pmod{\mathcal{L}_{\mathbf{k}}} \qquad  (\star)$$
for some $b_j\in R$ and stacked links $\mathbf{g}_j$ with $|\mathbf{g}_j|=|\mathbf{k}|$ and ${\rm cn}(\mathbf{g}_j)<{\rm cn}(\mathbf{k})$, such that each $\mathbf{g}_j$ is either a knot or has the form $\mathbf{k}_1\mathbf{k}_2$ with $\pi(\mathbf{k}_1)\cap\pi(\mathbf{k}_2)=\emptyset$.
\begin{enumerate}
  \item Suppose $\mathbf{k}$ contains a {\it convenient} arc, by which we mean $\mathbf{a}\in{\rm Ar}(\mathbf{k})$ such that
        ${\rm Cr}(\langle\mathbf{k}|\mathbf{a}\rangle,\mathbf{a})=\emptyset$ and ${\rm cn}(\mathbf{a})=1$, say
        ${\rm Cr}(\mathbf{a})=\{\mathsf{c}\}$.
        Then one of $\mathbf{k}^{\mathsf{c}}_\infty,\mathbf{k}^{\mathsf{c}}_0$ is a knot, and the other is congruent to a 2-component stacked link.
        For the EST $(\mathbf{k},\mathbf{k}^{\mathsf{c}}_\infty,\mathbf{k}^{\mathsf{c}}_0)$, by Lemma \ref{lem:EST},
        $$\varepsilon(\mathbf{k})\equiv q^{\frac{1}{2}}\varepsilon(\mathbf{k}^{\mathsf{c}}_\infty)
        +\overline{q}^{\frac{1}{2}}\varepsilon(\mathbf{k}^{\mathsf{c}}_0)\pmod{\mathcal{L}_{\mathbf{k}}}.$$
        Denote the 2-component stacked link in $\{\mathbf{k}^{\mathsf{c}}_\infty,\mathbf{k}^{\mathsf{c}}_0\}$ by $\mathbf{k}_1\mathbf{k}_2$. If
        $\pi(\mathbf{k}_1)\cap\pi(\mathbf{k}_2)\ne\emptyset$, then taking $\mathsf{c}'\in{\rm Cr}(\mathbf{k}_1,\mathbf{k}_2)$ to construct an EST,
        we obtain
        $$\varepsilon(\mathbf{k}_1\mathbf{k}_2)\equiv q^{\frac{1}{2}}\varepsilon(\mathbf{j}_1)
        +\overline{q}^{\frac{1}{2}}\varepsilon(\mathbf{j}_2)\pmod{\mathcal{L}_{\mathbf{k}}},$$
        where $\mathbf{j}_1=(\mathbf{k}_1\mathbf{k}_2)^{\mathsf{c}'}_\infty$, $\mathbf{j}_2=(\mathbf{k}_1\mathbf{k}_2)^{\mathsf{c}'}_0$,
        both being knots.
        Thus, $(\star)$ holds.
  \item In general, convenient arcs may not exist. However, we can always take $\mathbf{a}\in{\rm Ar}(\mathbf{k})$ with ${\rm cn}(\mathbf{a})=1$,
        then $\mathbf{k}^{\mathbf{a}}$ has $\mathbf{a}^\diamond$ as a convenient arc. So $(\star)$ holds for the knot $\mathbf{k}^{\mathbf{a}}$.
        On the other hand, by Lemma \ref{lem:crossing-change},
        $$\varepsilon(\mathbf{k})\equiv q^{\hat{\epsilon}(\mathbf{k},\mathbf{a})}\varepsilon(\mathbf{k}^{\mathbf{a}})
        +\varepsilon(\mathfrak{r}_u(\mathbf{k},\mathbf{a}))\pmod{\mathcal{L}_{\mathbf{k}}}.$$
        Thus, $(\star)$ holds for $\mathbf{k}$.
\end{enumerate}

Go on to deal with each $\mathbf{g}_j$. When $\mathbf{g}_j$ has the form $\mathbf{k}_1\mathbf{k}_2$ with $\pi(\mathbf{k}_1)\cap\pi(\mathbf{k}_2)=\emptyset$, we treat $\mathbf{k}_1,\mathbf{k}_2$ separately; observe that the procedure applied to $\mathbf{k}_1$ does not interfere that applied to $\mathbf{k}_2$.

Repeat such procedures. Ultimately we obtain $\varepsilon(\mathbf{k})\equiv\varepsilon(g)\pmod{\mathcal{L}_{\mathbf{k}}}$, where $g$ is a $R$-linear combination of multicurves of degree $|\mathbf{k}|$. Since $g=\mathbf{k}=\Theta(\mathbf{k})$ in $\mathcal{S}_n$, by Remark \ref{rmk:multicurve}, $\varepsilon(g)\equiv\varepsilon(\Theta(\mathbf{k}))\pmod{\mathcal{L}_{|\mathbf{k}|,1}}$.
Thus, $\varepsilon(\mathbf{k})\equiv\varepsilon(\Theta(\mathbf{k}))\pmod{\mathcal{L}_{\mathbf{k}}}$.
\end{proof}

\begin{figure}[H]
  \centering
  \includegraphics[width=12.8cm]{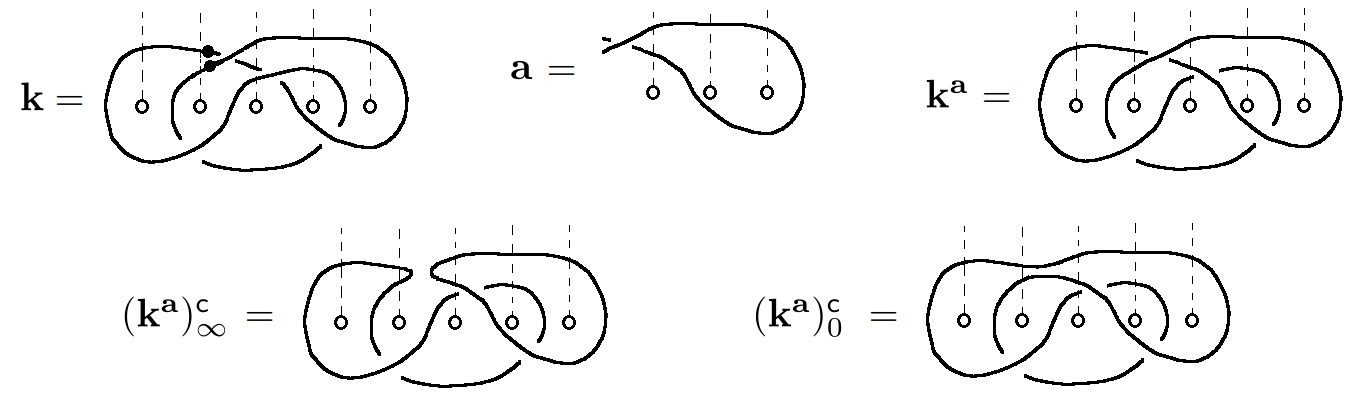}\\
  \caption{First row: pull $\mathbf{a}\in{\rm Ar}(\mathbf{k})$ up to the top, then $\mathbf{a}^\diamond\in{\rm Ar}(\mathbf{k}^{\mathbf{a}})$ is convenient.
  Second row: $(\mathbf{k}^{\mathbf{a}})^{\mathsf{c}}_\infty$ and $(\mathbf{k}^{\mathbf{a}})^{\mathsf{c}}_0$, obtained by resolving the unique crossing $\mathsf{c}\in{\rm Cr}(\mathbf{a}^\diamond)$.}\label{fig:convenient}
\end{figure}

\begin{exmp}
\rm The knot $\mathbf{k}$ given at the upper-left corner of Figure \ref{fig:convenient} admits no convenient arc, but it contains an arc $\mathbf{a}$ (the one bounded by the bullets) with the property that $\pi(\mathbf{a})$ has exactly one self-intersection. Resulted from pulling $\mathbf{a}$ up to the top, $\mathbf{a}^\diamond$ is a convenient arc of $\mathbf{k}^{\mathbf{a}}$.
Let $\mathsf{c}$ denote the unique crossing of $\mathbf{a}^\diamond$.
Vertically isotope $(\mathbf{k}^{\mathbf{a}})^{\mathsf{c}}_\infty$ into a stacked link $\mathbf{k}_1\mathbf{k}_2$, then
$(\mathbf{k}^{\mathbf{a}},\mathbf{k}_1\mathbf{k}_2,(\mathbf{k}^{\mathbf{a}})^{\mathsf{c}}_0)$ is an EST.
\end{exmp}

\begin{lem}\label{lem:simplify}
$\varepsilon(\mathcal{L}_{m,c+1})\subseteq\mathcal{L}_{m,c}$ for any $c\ge 0$.
\end{lem}

\begin{proof}
For type (i) elements of $\mathcal{L}_{m,c+1}$, suppose $g=\sum_ia_i\mathbf{k}_i\in\tilde{\theta}^{-1}(0)$, where $a_i\in R$ and $\mathbf{k}_i$ is a knot with $\|\mathbf{k}_i\|\le(m,c)$. By Lemma \ref{lem:resolve},
$\varepsilon(\mathbf{k}_i)\equiv\varepsilon(\Theta(\mathbf{k}_{i}))\pmod{\mathcal{L}_{m,c}}$ for each $i$.
Hence $\varepsilon(g)\equiv\varepsilon(\Theta(g))\equiv 0\pmod{\mathcal{L}_{m,c}}$.

For type (ii) elements, say
$$g=\mathbf{k}-{\sum}_sa_s(\mathbf{k}^{\mathbf{a}}|\mathbf{a}^\diamond|\mathbf{c}_s)-\mathfrak{r}_u(\mathbf{k},\mathbf{a}),$$
with $\|\mathbf{k}\|\preceq(m,c)$ and ${\rm ch}_u(\mathbf{a}^\diamond)=\sum_sa_s[\mathbf{c}_s]$.
Let $\mathbf{r}_s=(\mathbf{k}^{\mathbf{a}}|\mathbf{a}^\diamond|\mathbf{c}_s)$.
Then $a_s(\mathbf{r}_s-\varepsilon(\mathbf{r}_s))$ is a type (iv) element of $\mathcal{L}_{m,c}$ for each $s$.
By Remark \ref{rmk:chop-up} (ii), $\mathbf{j}\equiv\varepsilon(\mathbf{j})\pmod{\mathcal{L}_{m,c}}$ for each knot $\mathbf{j}$ appearing in $\mathfrak{r}_u(\mathbf{k},\mathbf{a})$.
Hence
\begin{align*}
\varepsilon(\mathbf{k})&\equiv\varepsilon_u(\mathbf{k},\mathbf{a})\equiv{\sum}_sa_s\mathbf{r}_s+\mathfrak{r}_u(\mathbf{k},\mathbf{a})  \\
&\equiv{\sum}_sa_s\varepsilon(\mathbf{r}_s)+\varepsilon(\mathfrak{r}_u(\mathbf{k},\mathbf{a}))
=\varepsilon\left({\sum}_sa_s\mathbf{r}_s+\mathfrak{r}_u(\mathbf{k},\mathbf{a})\right)\pmod{\mathcal{L}_{m,c}},
\end{align*}
where the first and second $\equiv$'s hold by definition. Thus, $\varepsilon(g)\in\mathcal{L}_{m,c}$.

Similarly for type (iii) elements.

Finally, let $\mathcal{L}'_{m,c}$ denote the submodule of $\mathcal{L}_{m,c}$ generated by type (iv) elements. Then clearly,
$\varepsilon(\mathcal{L}'_{m,c+1})\subseteq\mathcal{L}'_{m,c+1}=\mathcal{L}'_{m,c}$.
\end{proof}

\begin{nota}
\rm For each stacked link $\mathbf{l}$, there exists $N$ such that $\varepsilon^{N+k}(\mathbf{l})\equiv\varepsilon^N(\mathbf{l})\pmod{\mathcal{L}_{\mathbf{l}}}$ for all $k>0$; let $\varepsilon^\infty(\mathbf{l})\equiv\varepsilon^N(\mathbf{l})\pmod{\mathcal{L}_{\mathbf{l}}}$.
\end{nota}

\begin{lem}\label{lem:last}
$\varepsilon^\infty(\mathbf{l})\equiv\varepsilon^\infty(\Theta(\mathbf{l}))\pmod{\mathcal{L}_{|\mathbf{l}|,0}}$ for each stacked link $\mathbf{l}$.
\end{lem}

\begin{proof}
Suppose $\mathbf{l}=\mathbf{k}_1\cdots\mathbf{k}_r$, with each $\mathbf{k}_i$ a knot. We use induction on $r$ to prove the assertion.

If $r=1$, i.e. $\mathbf{l}$ is a knot, then by Lemma \ref{lem:resolve}, Remark \ref{rmk:multicurve} and Lemma \ref{lem:simplify}, $\varepsilon^\infty(\mathbf{l})\equiv\varepsilon^\infty(\Theta(\mathbf{l}))\pmod{\mathcal{L}_{|\mathbf{l}|,0}}$.

Suppose $r=\ell\ge 2$ and that the assertion holds when $r=\ell-1$.

If $\pi(\mathbf{k}_i)\cap\pi(\mathbf{k}_j)=\emptyset$ for all $i,j$,
then
$\Theta(\mathbf{l})=\Theta(\mathbf{k}_1)\cdots\Theta(\mathbf{k}_r)$ in $\mathcal{S}_n$.
By Lemma \ref{lem:resolve} and Remark \ref{rmk:multicurve},
$$\varepsilon(\mathbf{l})\equiv\varepsilon(\mathbf{k}_1)\cdots\varepsilon(\mathbf{k}_r)\equiv
\varepsilon(\Theta(\mathbf{k}_1))\cdots\varepsilon(\Theta(\mathbf{k}_r))\equiv\varepsilon(\Theta(\mathbf{l}))\pmod{\mathcal{L}_{|\mathbf{l}|,1}}.$$
Hence $\varepsilon^\infty(\mathbf{l})\equiv\varepsilon^\infty(\Theta(\mathbf{l}))\pmod{\mathcal{L}_{|\mathbf{l}|,0}}$.

From now on, assume $\pi(\mathbf{k}_i)\cap\pi(\mathbf{k}_j)\ne\emptyset$ for some $i\ne j$.
By Lemma \ref{lem:disjoint}, $\varepsilon(\mathbf{k}_s)\equiv\varepsilon(\mathbf{k}_t)\pmod{\mathcal{L}_{|\mathbf{l}|,0}}$ if $\pi(\mathbf{k}_s)\cap\pi(\mathbf{k}_t)=\emptyset$. Hence interchanging adjacent knots if necessary,
we may just assume $\pi(\mathbf{k}_1)\cap\pi(\mathbf{k}_2)\ne\emptyset$. Take $\mathsf{c}\in{\rm Cr}(\mathbf{k}_1,\mathbf{k}_2)$ to construct an EST $(\mathbf{g},\mathbf{g}^{\mathsf{c}}_\infty,\mathbf{g}^{\mathsf{c}}_0)$, with $\mathbf{g}=\mathbf{k}_1\mathbf{k}_2$, and
$\mathbf{g}^{\mathsf{c}}_\infty,\mathbf{g}^{\mathsf{c}}_0$ being knots.
By Lemma \ref{lem:EST},
$$\varepsilon(\mathbf{g})\equiv q^{\frac{1}{2}}\varepsilon(\mathbf{g}^{\mathsf{c}}_\infty)
+\overline{q}^{\frac{1}{2}}\varepsilon(\mathbf{g}^{\mathsf{c}}_0)\pmod{\mathcal{L}_{\mathbf{g}}}.$$
Let $\mathbf{l}_1=\mathbf{g}^{\mathsf{c}}_\infty\mathbf{k}_3\cdots\mathbf{k}_r$, $\mathbf{l}_2=\mathbf{g}^{\mathsf{c}}_0\mathbf{k}_3\cdots\mathbf{k}_r$, each having $\ell-1$ components.
Then
$$\varepsilon(\mathbf{l})\equiv q^{\frac{1}{2}}\varepsilon(\mathbf{l}_1)+\overline{q}^{\frac{1}{2}}\varepsilon(\mathbf{l}_2)\pmod{\mathcal{L}_{|\mathbf{l}|,0}}.$$
By the inductive hypothesis, $\varepsilon^\infty(\mathbf{l}_i)\equiv\varepsilon^\infty(\Theta(\mathbf{l}_i))\pmod{\mathcal{L}_{|\mathbf{l}|,0}}$.
Hence
$$\varepsilon^\infty(\mathbf{l})\equiv q^{\frac{1}{2}}\varepsilon^\infty(\Theta(\mathbf{l}_1))+\overline{q}^{\frac{1}{2}}\varepsilon^\infty(\Theta(\mathbf{l}_2))
\equiv\varepsilon^\infty(\Theta(\mathbf{l}))\pmod{\mathcal{L}_{|\mathbf{l}|,0}}.$$
The second $\equiv$ is due to that $\Theta(\mathbf{l})=q^{\frac{1}{2}}\Theta(\mathbf{l}_1)+\overline{q}^{\frac{1}{2}}\Theta(\mathbf{l}_2)$ in $\mathcal{S}_n$; we refer to Remark \ref{rmk:multicurve}.
\end{proof}

\section{Proof of Theorem \ref{thm:main}}

\begin{defn}
\rm By a {\it monomial} we mean a product of elements of $G$, which can be regarded as a stacked link.
\end{defn}

Let $\mathcal{J}$ denote the two-sided ideal of $\mathcal{T}$ generated by elements of the form $\sum_ia_i\mathbf{u}_i$ such that $a_i\in R$, and $\mathbf{u}_i$ is a monomial with $|\mathbf{u}_i|\le 2h$, and $\sum_ia_i[\mathbf{u}_i]=0$.

\begin{lem}\label{lem:J}
$\varepsilon^\infty(\mathcal{L}_{m,0})\subset\mathcal{J}$ for each $m$.
\end{lem}

\begin{proof}
We use induction on $m$ to prove the assertion, which obviously holds when $m\le 2h+1$.
Suppose $m>2h+1$ and $\varepsilon^\infty(\mathcal{L}_{m',0})\subset\mathcal{J}$ for $m'<m$.

Recall that $\mathcal{L}_{m,0}$ is generated by elements of the form $\mathbf{j}(\sum_ia_i\mathbf{l}_i)\mathbf{j}'$ such that $\sum_ia_i[\mathbf{l}_i]=0$ and each $\mathbf{l}_i$ is a stacked link with $|\mathbf{l}_i|<m$.

For such $\sum_ia_i\mathbf{l}_i$ (satisfying $\sum_ia_i\Theta(\mathbf{l}_i)=0$ in $\mathcal{V}$), let $m'=\max_i\{|\mathbf{l}_i|\}<m$, then by Lemma \ref{lem:last}, $\varepsilon^\infty(\mathbf{l}_i)\equiv\varepsilon^\infty(\Theta(\mathbf{l}_i))\pmod{\mathcal{L}_{m',0}}$.
Hence
$$\varepsilon^\infty\left({\sum}_ia_i\mathbf{l}_i\right)\equiv{\sum}_ia_i\varepsilon^\infty(\Theta(\mathbf{l}_i))
=\varepsilon^\infty\left({\sum}_ia_i\Theta(\mathbf{l}_i)\right)\equiv 0\pmod{\mathcal{L}_{m',0}}.$$
Consequently,
$$\varepsilon^\infty\left({\sum}_ia_i\mathbf{l}_i\right)=\varepsilon^\infty\left(\varepsilon^\infty\left({\sum}_ia_i\mathbf{l}_i\right)\right)
\in\varepsilon^\infty(\mathcal{L}_{m',0})\subset\mathcal{J}.$$
Thus, $\varepsilon^\infty(\mathbf{j}(\sum_ia_i\mathbf{l}_i)\mathbf{j}')\in\mathcal{J}$.
\end{proof}


\begin{thm}\label{thm:bound}
The skein algebra $\mathcal{S}_n$ is generated by $G$, and the ideal of defining relations of $\mathcal{S}_n$ is $\mathcal{J}$, i.e., $\ker\theta=\mathcal{J}$.
\end{thm}

\begin{proof}
We have $\theta(\varepsilon^\infty(\mathbf{m}))=[\mathbf{m}]$ for each multicurve $\mathbf{m}$.
Since $\mathcal{S}_n$ is spanned by multicurves, $\theta$ is surjective, i.e. $\mathcal{S}_n$ is generated by $G$.

Suppose $\theta(f)=0$ for $f=\sum_i\beta_i\mathbf{u}_i\in\mathcal{T}$, with $\beta_i\in R$ and $\mathbf{u}_i$ a monomial.
Then $f\in\mathcal{L}_{m,0}$ for $m=1+\max_i\{|\mathbf{u}_i|\}$, so by Lemma \ref{lem:J}, $\varepsilon^\infty(f)\in\mathcal{J}$.
On the other hand, $\varepsilon^\infty(f)=f$, hence $f\in\mathcal{J}$. Thus, $\ker\theta=\mathcal{J}$.
\end{proof}

\begin{nota}
\rm Given $\mathsf{w}=(w_1,\ldots,w_n)\in(\mathbb{Z}_{\ge0})^n$, let $|\mathsf{w}|=\sum_{i=1}^nw_i$, and let $\mathcal{Q}_{\mathsf{w}}$ denote the $R$-submodule of $\mathcal{T}$ generated by monomials $\mathbf{v}$ with $|\mathbf{v}|_i\le w_i$ for each $i$.
Let $\Delta_{\mathsf{w}}=\Delta_{s_1,\ldots,s_k}$ if $s_1<\cdots<s_k$ are the $i$'s with $w_i>0$.
\end{nota}
Note that elements of $\mathcal{Q}_{\mathsf{w}}$ represent elements of $\mathcal{S}(\Delta_{\mathsf{w}};R)$.

Let $\mathcal{I}$ denote the two-sided ideal of $\mathcal{T}$ generated by $\bigcup_{|\mathsf{w}|\le 2h}\mathcal{C}_{\mathsf{w}}$, where
$$\mathcal{C}_{\mathsf{w}}=\big\{f\in\mathcal{Q}_{\mathsf{w}}\colon f=0\ \text{in\ }\mathcal{S}(\Delta_{\mathsf{w}};R)\big\}.$$

Given $f\in\mathcal{C}_{\mathsf{w}}$ with $|\mathsf{w}|\le 2h$, say $f={\sum}_{i=1}^r\beta_i\mathbf{u}_i$, with $\beta_i\in R$ and $\mathbf{u}_i$ a monomial, by definition $|\mathbf{u}_i|\le|\mathsf{w}|\le 2h$ for all $i$.
In this sense, we say that $f=0$ is a relation of degree $\le 2h$ {\it supported} by the subsurface $\Delta_{\mathsf{w}}$.

\begin{lem}\label{lem:ideal}
$\mathcal{J}=\mathcal{I}$.
\end{lem}

\begin{proof}
By definition, $\mathcal{I}\subseteq\mathcal{J}$. It suffices to show $\mathcal{J}\subseteq\mathcal{I}$.

For each stacked link $\mathbf{l}$, let $\mathsf{w}=(|\mathbf{l}|_1,\ldots,|\mathbf{l}|_n)$, then $\varepsilon^\infty(\mathbf{l})$ and $\varepsilon^\infty(\Theta(\mathbf{l}))$ can be defined totally within $\Delta_{\mathsf{w}}$, and $\varepsilon^\infty(\mathbf{l})\equiv\varepsilon^\infty(\Theta(\mathbf{l}))\pmod{\mathcal{C}_{\mathsf{w}}}$.
When $|\mathbf{l}|\le 2h$, we have $\mathcal{C}_{\mathsf{w}}\subset\mathcal{I}$, so $\varepsilon^\infty(\mathbf{l})\equiv\varepsilon^\infty(\Theta(\mathbf{l}))\pmod{\mathcal{I}}$.

If $f\in\mathcal{J}$, say $f=\sum_ia_i\mathbf{u}_i\in\theta^{-1}(0)$ with $|\mathbf{u}_i|\le 2h$, then $\Theta(f)=0$, i.e. the irreducible multicurves obtained from the $\mathbf{u}_i$'s cancel out, hence
$$f=\varepsilon^\infty(f)\equiv\varepsilon^\infty(\Theta(f))\equiv 0\pmod{\mathcal{I}}.$$
\end{proof}

When $\alpha^{-1}\in R$, so that $h=3$, Theorem \ref{thm:bound} together with Lemma \ref{lem:ideal} establish Theorem \ref{thm:main} (a).

When $\alpha^{-1}\notin R$, so that $h=n+1$, the ideal $\mathcal{I}$ can be described more sharply.
For each $\mathsf{w}=(w_1,\ldots,w_n)$, if $s_1<\cdots<s_k$ are the $i$'s with $w_i>0$, then under the diffeomorphism $\Delta_{\mathsf{w}}\to\Sigma_{0,k+1}$, $\mathcal{S}(\Delta_{\mathsf{w}};R)\cong\mathcal{S}_k$, and $\mathcal{C}_{\mathsf{w}}$ is identified with a submodule of the defining ideal $\mathcal{I}_k$ of $\mathcal{S}_k$. Applying Theorem \ref{thm:bound} and Lemma \ref{lem:ideal} to $n=k$, we see that $\mathcal{I}_k$ is generated by relations of degree $\le 2k+2$, which, viewed as relations in $\mathcal{S}_n$, are supported by $\Delta_{s_1,\ldots,s_k}$.
This completes the proof of Theorem \ref{thm:main} (b).

\section{Appendix}

\subsection{Proof of Lemma \ref{lem:reduce-multi-curve}}\label{sec:irreducible}

\begin{proof}
Suppose $\phi=\{\phi_t\}_{t\in[0,1]}$ is an isotopy with $\phi_1(\mathbf{m}_0)=\mathbf{m}_1$. We show that it can be modified into a fine isotopy.

Let $Y_k=\mathbf{z}_k\times[0,1]$, and let $Y=\cup_{k=1}^nY_k$.
The ``worldsheet" of $\phi$,
$$\Omega(\phi):=\{(\phi_t(\mathsf{x}),t)\colon \mathsf{x}\in\mathbf{m}_0, t\in[0,1]\}$$
is a disjoint union of annuli embedded in $\Sigma\times[0,1]$.

\begin{figure}[h]
  \centering
  \includegraphics[width=10.5cm]{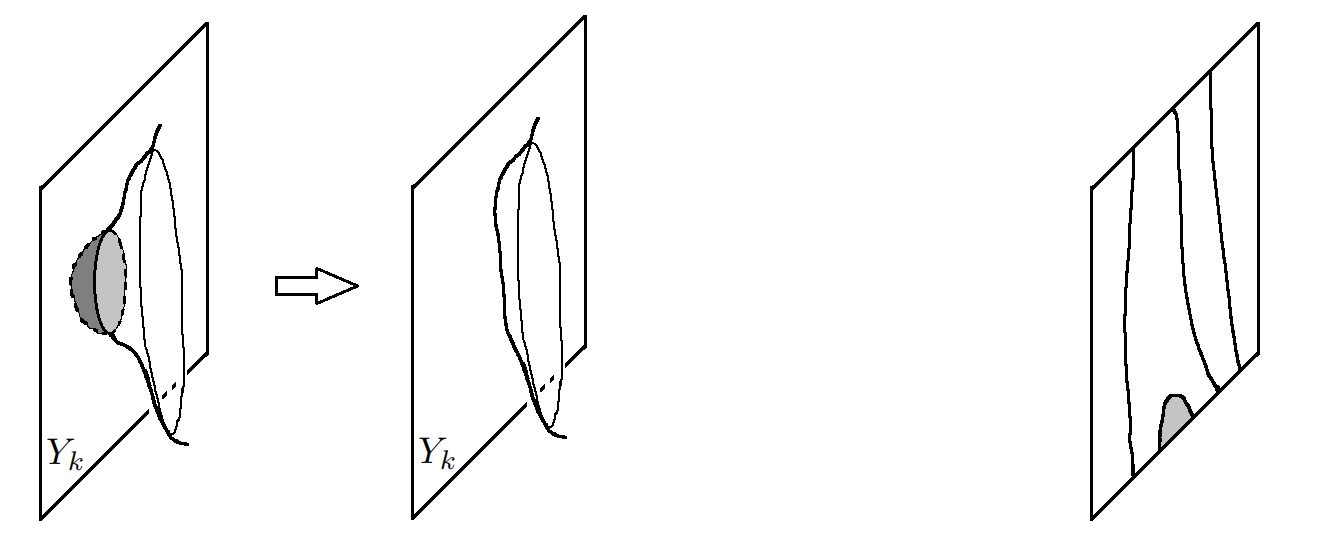}\\
  \caption{Left: the shaded region is the ball $B$ bounded by $D\cup_CE$; squeezing $B$ out, we eliminate a circle in $\Omega\cap Y$.
  Right: under the continuous map $(\mathsf{x},t)\mapsto\omega_t^{-1}(\mathsf{x})$, the cap becomes an arc $\mathbf{c}\subset\mathbf{m}_0$, and the image of the shaded disk provides an endpoint-fixing homotopy between $\mathbf{c}$ and a segment of $\mathbf{z}_k$.}\label{fig:squeeze}
\end{figure}

If $\Omega(\phi)\cap Y$ contains at least one circle, then each circle bounds a disk in $\Omega(\phi)$. Let $D$ be such a disk which is minimal with respect to inclusion.
Suppose $C=\partial D\subset\Omega\cap Y_k$ for some $k$, then $C$ also bounds some disk $E$ in $Y_k$. Since $D\cup_CE$ is a 2-sphere, by
the {\it generalized Schoenflies Theorem} it bounds a ball $B\subset\Sigma\times(0,1)$. Due to the minimality of $D$, the interior of $B$ does not intersect with $\Omega(\phi)$.
We can horizontally squeeze $B$ out, so as to push $D$ away from $Y_k$, as illustrated in the left part of Figure \ref{fig:squeeze}. The effect is equivalent to replacing $\phi$ by another isotopy $\psi=\{\psi_t\}_{t\in[0,1]}$, such that
$\Omega(\psi)\cap Y$ contains less circles than $\Omega(\phi)\cap Y$.

Repeating this procedure, ultimately we obtain an isotopy $\omega=\{\omega_t\}_{t\in[0,1]}$ such that $\omega_1(\mathbf{m}_0)=\mathbf{m}_1$ and $\Omega(\omega)\cap Y$ contains no circle.

For each $k$, we claim that $\Omega(\omega)\cap Y_k$ does not contain a cup or cap. Indeed, if there was a cap (as illustrated in the right part of Figure \ref{fig:squeeze}), then $\mathbf{m}_0$ would contain an arc homotopic to a segment of $\mathbf{z}_k$, where the homotopy fixes endpoints. This contradicts the irreducibility of $\mathbf{m}_0$.
Similarly, the irreducibility of $\mathbf{m}_1$ forbids $\Omega(\omega)\cap Y_k$ to contain a cup.

Thus, $\Omega(\omega)$ is a surface diffeomorphic to a disjoint union of annuli such that $\Omega(\omega)\cap Y$ provides a trivial cobordism between $\mathbf{m}_0\cap\mathbf{z}$ and $\mathbf{m}_1\cap\mathbf{z}$.
Reparameterizing $\Omega(\omega)$, we can obtain a fine isotopy $\varphi=\{\varphi_t\}_{t\in[0,1]}$ with $\Omega(\varphi)=\Omega(\omega)$; in particular, $\varphi_1(\mathbf{m}_0)=\mathbf{m}_1$.
\end{proof}

\subsection{Proof of Lemma \ref{lem:substitution}}

\begin{proof}
We prove (i) by constructing ${\rm ch}_u(\mathbf{a})$. The proof for (ii) is parallel.

When $\mathbf{a}$ is simple (so that it can be identified with an arc in $\Sigma$) and contains a shrinkable subarc $\mathbf{b}$, just
let ${\rm ch}_u(\mathbf{a})\in\mathcal{S}(\mathsf{a},\mathsf{b})$ be the element represented by the arc obtained by reducing $\mathbf{b}$.

When $\mathbf{a}$ is non-simple, we can resolve its crossings and reduce shrinkable subarcs to obtain $\sum_rc_r\mathbf{a}_r$, for some $c_r\in\mathcal{T}$ and simple irreducible arcs $\mathbf{a}_r$. As is easy to see, $|c_r|_i+|\mathbf{a}_r|_i\le|\mathbf{a}|_i$ for each $i$.
Put ${\rm ch}_u(\mathbf{a})=\sum_rc_r{\rm ch}_u(\mathbf{a}_r)$, where ${\rm ch}_u(\mathbf{a}_r)=[\mathbf{a}_r]$ if $|\mathbf{a}_r|\le 2$, and ${\rm ch}_u(\mathbf{a}_r)$ will be defined in a moment if $|\mathbf{a}_r|=3$ (in which case $c_r\in R$).

From now on, assume $\mathbf{a}$ to be simple and irreducible, so that $\#{\rm supp}(\mathbf{a})\in\{2,3\}$.

If $\#{\rm supp}(\mathbf{a})=2$, say ${\rm supp}(\mathbf{a})=\{i,j\}$ with $i<j$, then either $|\mathbf{a}|_i=2$, $|\mathbf{a}|_j=1$, or
$|\mathbf{a}|_i=1$, $|\mathbf{a}|_j=2$. The two possibilities are respectively shown as the left-hand-side of the first and second row in Figure \ref{fig:basic-0}.
In either case, correspondingly let ${\rm ch}_u(\mathbf{a})$ be the expression given by the right-hand-side.

\begin{figure}[H]
  \centering
  \includegraphics[width=12cm]{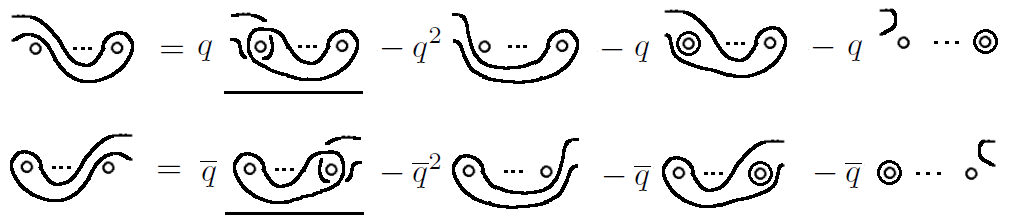}
  \caption{When ${\rm supp}(\mathbf{a})=\{i,j\}$ with $i<j$, a portion of $\mathbf{a}$ near $\mathsf{p}_{i},\mathsf{p}_{j}$ takes the form in the left-hand-side of either row. Each equation is easy to verify by resolving the crossings in the underlined term.}\label{fig:basic-0}
\end{figure}

Now suppose $\#{\rm supp}(\mathbf{a})=3$.

Let $\Sigma^{\rm cut}$ denote the surface obtained by cutting $\Sigma$ along $\mathbf{z}$. Fix a diffeomorphism $\Sigma^{\rm cut}\to D^2$ to identify $\Sigma^{\rm cut}$ with $D^2$. Let ${\rm gl}:D^2\to\Sigma$ denote the gluing map.
For $\mathbf{x}\subset\Sigma$, denote $\tilde{\mathbf{x}}$ for ${\rm gl}^{-1}(\mathbf{x})$.

Up to orientation-preserving diffeomorphism of $D^2$ to itself, $\tilde{\mathbf{a}}$ has the form of one of $\tilde{\mathbf{r}}_1,\ldots,\tilde{\mathbf{r}}_4$ (shown in Figure \ref{fig:piece}).
Indeed, $\mathbf{a}$ is divided by $\mathbf{z}$ into three parts, say $\mathbf{a}=\mathbf{c}_1\mathbf{b}\mathbf{c}_2$, with $\mathbf{c}_i\cap\mathbf{z}=\partial\mathbf{c}_i\cap\mathbf{z}$, $i=1,2$; depending on the relative positions of $\partial\widetilde{\mathbf{c}_1}$,
$\partial\widetilde{\mathbf{c}_2}$ and $\partial\tilde{\mathbf{b}}$, there are exactly four possibilities.

\begin{figure}[H]
  \centering
  \includegraphics[width=10.3cm]{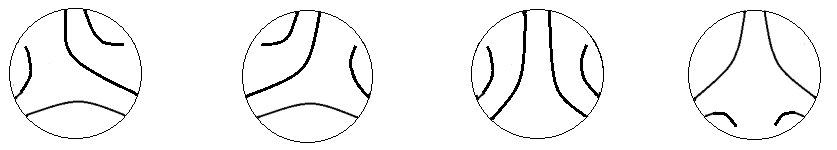}
  \caption{From left to right: $\tilde{\mathbf{r}}_1$, $\tilde{\mathbf{r}}_2$, $\tilde{\mathbf{r}}_3$, $\tilde{\mathbf{r}}_4$.}\label{fig:piece}
\end{figure}

\begin{figure}[H]
  \centering
  \includegraphics[width=11.3cm]{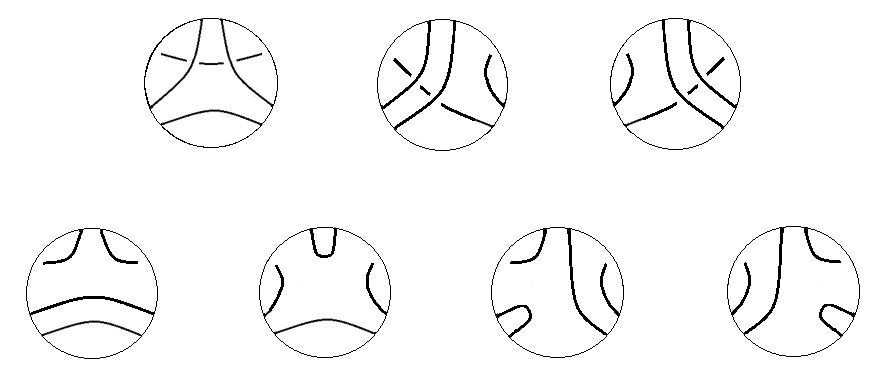}
  \caption{From left to right, first row: $\tilde{\mathbf{u}}_1$, $\tilde{\mathbf{u}}_2$, $\tilde{\mathbf{u}}_3$;
  second row: $\tilde{\mathbf{v}}_1$, $\tilde{\mathbf{v}}_2$, $\tilde{\mathbf{v}}_3$, $\tilde{\mathbf{v}}_4$.}\label{fig:term}
\end{figure}

\begin{figure}[H]
  \centering
  \includegraphics[width=13cm]{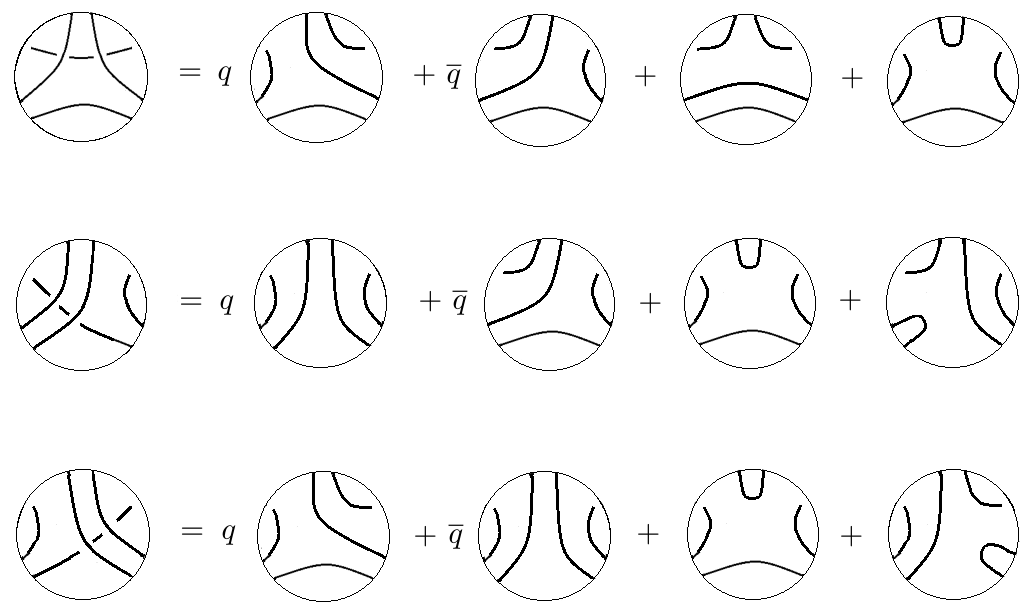}
  \caption{These are deduced by resolving crossings via skein relations.}\label{fig:basic-1}
\end{figure}

\begin{figure}[H]
  \centering
  \includegraphics[width=13cm]{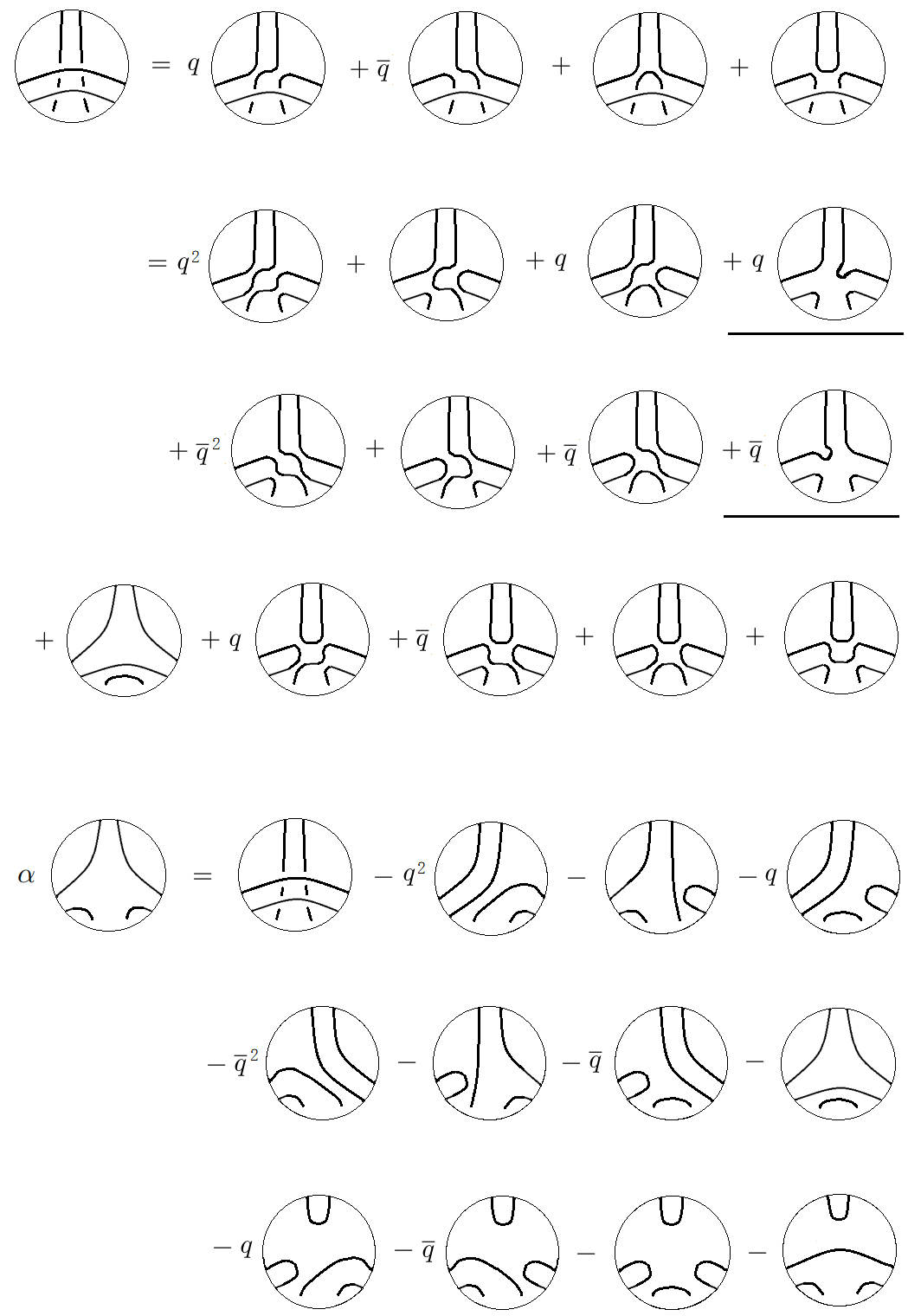}
  \caption{The upper equation is deduced by successively resolving crossings; moving all terms except the underlined ones to the left-hand-side yields the lower equation. When glued back to an equation in $\mathcal{S}(\mathsf{a},\mathsf{b})$, the first term is $\alpha[\mathbf{r}_4]$, and each of the other terms has the form $a[\mathbf{c}]$ with $a\in\mathcal{T}$ and $\mathbf{c}\in P(\mathsf{a},\mathsf{b})$.}\label{fig:basic-2}
\end{figure}

Introduce more elements in Figure \ref{fig:term}. The equations in Figure \ref{fig:basic-1} read
$$\left(\begin{array}{ccc} \tilde{\mathbf{u}}_1 \\ \tilde{\mathbf{u}}_2 \\ \tilde{\mathbf{u}}_3 \end{array}\right)
=\left(\begin{array}{ccc} q & \overline{q} & 0  \\ 0 & \overline{q} & q \\ q & 0 & \overline{q} \end{array}\right)
\left(\begin{array}{ccc} \tilde{\mathbf{r}}_1 \\ \tilde{\mathbf{r}}_2 \\ \tilde{\mathbf{r}}_3 \end{array}\right)
+\left(\begin{array}{ccc} \tilde{\mathbf{v}}_1+\tilde{\mathbf{v}}_2 \\ \tilde{\mathbf{v}}_2+\tilde{\mathbf{v}}_3 \\ \tilde{\mathbf{v}}_2+\tilde{\mathbf{v}}_4
\end{array}\right),$$
which can be solved as
\begin{align}
\tilde{\mathbf{r}}_1&=\alpha^{-1}\big(\overline{q}^2\tilde{\mathbf{u}}_1-\overline{q}^2\tilde{\mathbf{u}}_2+\tilde{\mathbf{u}}_3-\overline{q}^2\tilde{\mathbf{v}}_1
-\tilde{\mathbf{v}}_2+\overline{q}^2\tilde{\mathbf{v}}_3-\tilde{\mathbf{v}}_4\big), \nonumber  \\
\tilde{\mathbf{r}}_2&=\alpha^{-1}\big(q^2\tilde{\mathbf{u}}_1+\tilde{\mathbf{u}}_2-q^2\tilde{\mathbf{u}}_3-q^2\tilde{\mathbf{v}}_1-\tilde{\mathbf{v}}_2
-\tilde{\mathbf{v}}_3+q^2\tilde{\mathbf{v}}_4 \big),  \nonumber  \\
\tilde{\mathbf{r}}_3&=\alpha^{-1}(-\tilde{\mathbf{u}}_1+\tilde{\mathbf{u}}_2+\tilde{\mathbf{u}}_3+\tilde{\mathbf{v}}_1-\tilde{\mathbf{v}}_2
-\tilde{\mathbf{v}}_3-\tilde{\mathbf{v}}_4).  \label{eq:r3}
\end{align}

Let $\mathbf{u}_i={\rm gl}(\tilde{\mathbf{u}}_i)$, $\mathbf{v}_i={\rm gl}(\tilde{\mathbf{v}}_i)$ and $\mathbf{r}_i={\rm gl}(\tilde{\mathbf{r}}_i)$.
Observe that each of $\mathbf{u}_i$, $\mathbf{v}_i$ has the form $a\mathbf{c}$, with $a\in\mathcal{T}$ and $\mathbf{c}\in P(\mathsf{a},\mathsf{b})$.
Thus, for each $1\le i\le 3$, there exist $a_s\in\mathcal{T}$ and $\mathbf{c}_s\in P(\mathsf{a},\mathsf{b})$ such that $[\mathbf{r}_i]=\sum_sa_s[\mathbf{c}_s]$ in $\mathcal{S}(\mathsf{a},\mathsf{b})$.

From Figure \ref{fig:basic-2} we see that $[\mathbf{r}_4]$ also equals such a combination; note that the image of the beginning term in Figure \ref{fig:basic-2} under ${\rm gl}$ already has the form $a\mathbf{c}$, with $a\in\mathcal{T}$, $\mathbf{c}\in P(\mathsf{a},\mathsf{b})$.

Finally, it is easy to verify that $|a_s|_i+|\mathbf{c}_s|_i\le |\mathbf{a}|_i$ for all $s,i$.
\end{proof}

\begin{figure}[H]
  \centering
  \includegraphics[width=12.5cm]{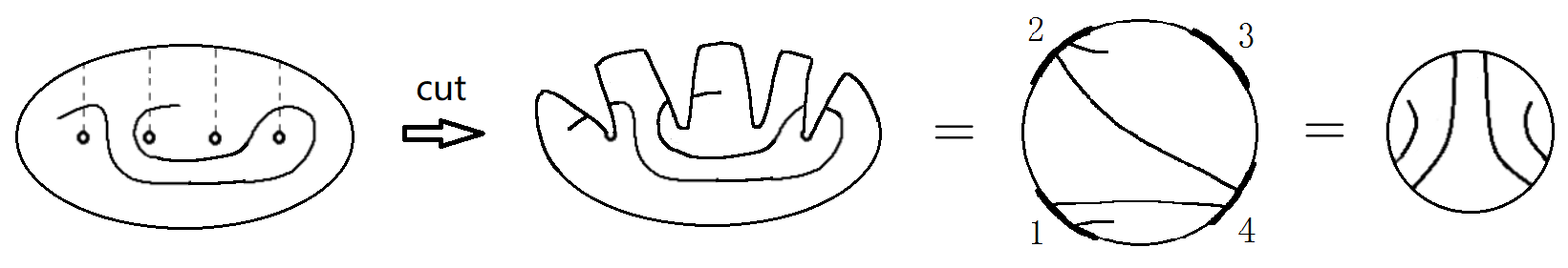}
  \caption{When $\Sigma$ is cut into $\Sigma^{\rm cut}\cong D^2$, a simple irreducible arc $\mathbf{a}$ with $|\mathbf{a}|=\#{\rm supp}(\mathbf{a})=3$ becomes a disjoint union of simple arcs in $D^2$. At the lower-right corner, the segment in bold labeled with $k$ stands for ${\rm gl}^{-1}(\mathbf{z}_{k})$.}\label{fig:symbol}
\end{figure}

\begin{figure}[H]
  \centering
  \includegraphics[width=13cm]{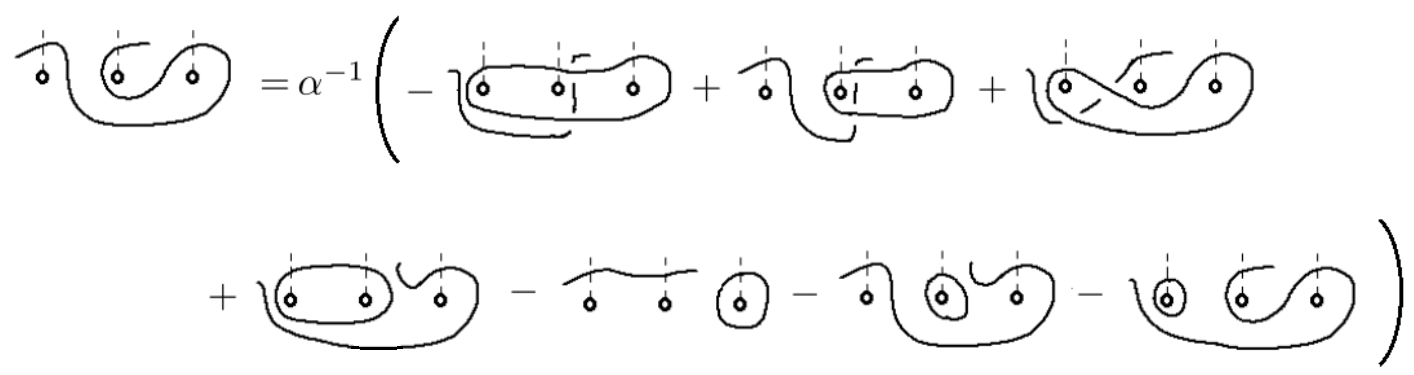}
  \caption{The arc $\mathbf{a}$ is chopped up into ${\rm ch}_u(\mathbf{a}^\diamond)$. Here $\mathbf{z}_3$ is omitted.}\label{fig:curve}
\end{figure}

\begin{exmp}
\rm Shown at the upper-left corner of Figure \ref{fig:symbol} is a simple irreducible arc $\mathbf{a}\subset\Sigma_{0,5}$ with
${\rm supp}(\mathbf{a})=\{1,2,4\}$. The lower-left corner has the form $\tilde{\mathbf{r}}_3$.
As Figure \ref{fig:curve} shows, $\mathbf{a}$ can be chopped up, in a form parallel to the formula (\ref{eq:r3}) for $\tilde{\mathbf{r}}_3$.
\end{exmp}

\begin{exmp}
\rm Fix $\mathsf{a},\mathsf{b}\in\Sigma\times\{0\}$.
Let $\mathbf{c}\in F_0(\mathsf{a},\mathsf{b})$ be an oriented simple arc with ${\rm word}(\mathbf{c})=x_{i_1}^{\nu_1}\cdots x_{i_m}^{\nu_m}$. When $m\ge 1$, we denote $\mathbf{c}$ as $\mathbf{x}_{i^\ast_1\cdots i^\ast_m}$, with $i_k^\ast=i_k$ (resp. $i_k^\ast=\overline{i_k}$ if $\nu_k=1$ (resp. $\nu_k=-1$); when $m=0$, denote $\mathbf{c}$ as $\mathbf{e}$.

Let $1\le i<j\le n$. Holding in $\mathcal{S}(\mathsf{a},\mathsf{b})$, the equations in Figure \ref{fig:basic-0} read
\begin{align*}
\mathbf{x}_{ij\overline{i}}&=q(t_{ij}\mathbf{x}_i-q\mathbf{x}_j-t_i\mathbf{x}_{ij}-t_j\mathbf{e}),  \\
\mathbf{x}_{\overline{j}ij}&=\overline{q}(t_{ij}\mathbf{x}_j-\overline{q}\mathbf{x}_i-t_j\mathbf{x}_{ij}-t_i\mathbf{e}).
\end{align*}

\begin{figure}[h]
  \centering
  \includegraphics[width=13cm]{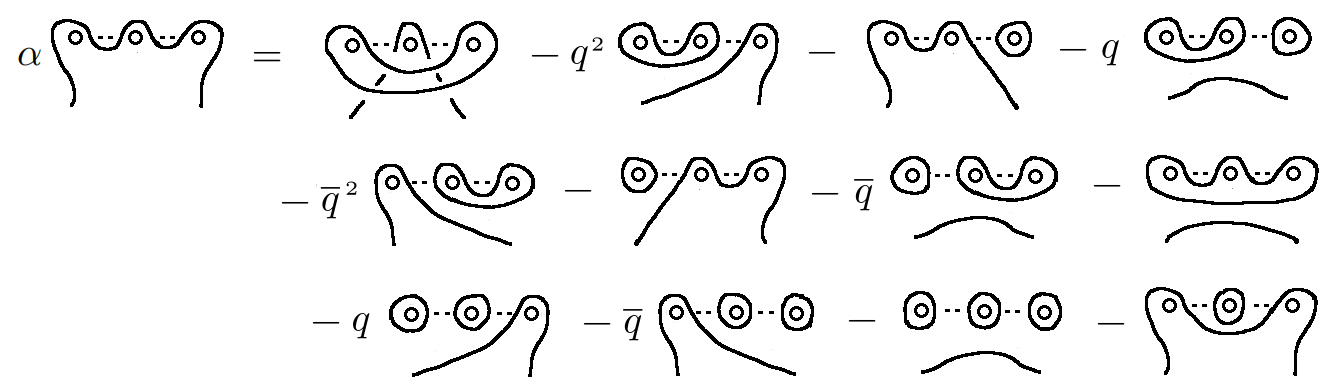}
  \caption{This is a possible case of the equation in Figure \ref{fig:basic-2}.}\label{fig:expand}
\end{figure}

Suppose $1\le i<j<k\le n$. The equation in Figure \ref{fig:expand} reads
\begin{align*}
\alpha\mathbf{x}_{ijk}=\ &t_{ik}\mathbf{x}_j-q^2t_{ij}\mathbf{x}_k-t_k\mathbf{x}_{ij}-qt_{ij}t_k\mathbf{e}-\overline{q}^2t_{jk}\mathbf{x}_i-t_i\mathbf{x}_{jk}
-\overline{q}t_it_{jk}\mathbf{e}-t_{ijk}\mathbf{e} \\
&-qt_it_j\mathbf{x}_k-\overline{q}t_jt_k\mathbf{x}_i-t_it_jt_k\mathbf{e}-t_j\mathbf{x}_{ik}.
\end{align*}
Hence the following holds in $\mathcal{S}(\mathsf{a},\mathsf{b})$:
\begin{align*}
\mathbf{x}_{ijk}=\alpha^{-1}\big(&t_{ik}\mathbf{x}_j-t_i\mathbf{x}_{jk}-t_j\mathbf{x}_{ik}-t_k\mathbf{x}_{ij}
-(\overline{q}^2t_{jk}+\overline{q}t_jt_k)\mathbf{x}_i  \\
&-(q^2t_{ij}+qt_it_j)\mathbf{x}_k-(t_{ijk}+qt_{ij}t_k+\overline{q}t_it_{jk}+t_it_jt_k)\mathbf{e}\big).
\end{align*}

\end{exmp}

\subsection{Proof of Lemma \ref{lem:shorten}}

\begin{figure}[h]
  \centering
  \includegraphics[width=12.5cm]{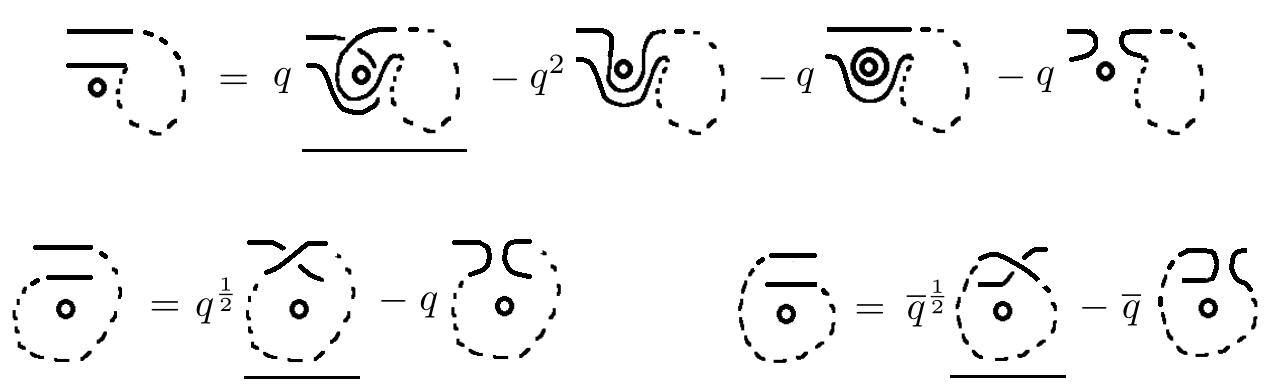}\\
  \caption{Each minimal shortenable arc $\mathbf{a}$ can be chopped up as $\sum_sa_s\mathbf{c}_s$, with $a_s\in\mathcal{T}$ and $\mathbf{c}_s$ unshortenable. These equations can be verified by resolving the crossings in the underlined terms.}\label{fig:shorten}
\end{figure}

\begin{proof}
We only prove (i); the proof for (ii) is similar.

Similarly as in the proof of roof of Lemma \ref{lem:substitution}, resolving the crossings of $\mathbf{a}$ if necessary, we may assume $\mathbf{a}$ to be simple.

Choose an orientation for $\mathbf{a}$.
Suppose ${\rm word}(\mathbf{a})=x_{i_1}^{\nu_1}\cdots x_{i_m}^{\nu_m}$, with $m\le n+1$ and $i_1=i_m$. The assertion is trivial if $m\le 2$.
Assume $m\ge 3$.

When $(\nu_1,\nu_m)=(1,-1)$, the assertion follows from the first row in Figure \ref{fig:shorten}.
Note that the underlined term has the form $\mathbf{j}\mathbf{c}$, where $\mathbf{c}\in Q(\mathsf{a},\mathsf{b})$ and $\mathbf{j}$ is a simple curve
with $|\mathbf{j}|_i\le 1$ for all $i$; such $\mathbf{j}$ belongs to $G$, so $\mathbf{j}\in\mathcal{T}$.

When $(\nu_1,\nu_m)=(1,1)$, there are two possibilities, and the assertion is clear from the second row in Figure \ref{fig:shorten}.

The cases when $(\nu_1,\nu_m)=(-1,-1)$ and $(\nu_1,\nu_m)=(-1,1)$ are similar.
\end{proof}





\bigskip

\noindent
Haimiao Chen (orcid: 0000-0001-8194-1264)\ \ \ \ {\it chenhm@math.pku.edu.cn} \\
Department of Mathematics, Beijing Technology and Business University, \\
Liangxiang Higher Education Park, Fangshan District, Beijing, China

\end{document}